\documentclass[11pt,a4paper, reqno]{amsart}  

\usepackage{latexsym}

\usepackage{amsmath,amsfonts,amssymb,amsthm,amscd,epsfig}  
\usepackage[english]{babel}  
\usepackage{indentfirst}  
\usepackage[mathscr]{eucal}  
\usepackage{graphics}  
  
\numberwithin{equation}{section}  
  
\usepackage{amsfonts,layout}
\newcommand{\N}{\mathbb N}
\newcommand{\R}{\mathbb R}
\def\QED{\begin{flushright}$\Box$\end{flushright}}

\def \lg{\langle} 
\def \rg{\rangle} 
\newcommand{\be}{\begin{equation}}
\newcommand{\ee}{\end{equation}}
\newcommand{\Sc}{{\mathbb S}}

\def\ep{\varepsilon}

\def\bx{{\bar x}}
\def\by{{\bar y}}

\newlength{\longarc}
\newcommand{\arc}[1]{\settowidth{\longarc}{$#1$}%
         \addtolength{\longarc}{-0.5em}%
\unitlength \longarc%
\ensuremath{%
         \stackrel{\begin{picture}(1,0.2)
                 \put(-5.5,-3){$\frown$}
         \end{picture}}
         {#1}\rule{0.01cm}{0cm}
}}
\newcommand{\conereg}{\arc{\mathcal C}}

\newlength{\longarcbis}

\newcommand{\cone}{{\widehat{\rule{2pt}{0pt}\mathcal{C}\rule{2pt}{0pt}}}}
\newtheorem{Theorem}{Theorem}[section]
\newtheorem{Definition}[Theorem]{Definition}
\newtheorem{Proposition}[Theorem]{Proposition}
\newtheorem{Lemma}[Theorem]{Lemma}
\newtheorem{Corollary}[Theorem]{Corollary}
\newtheorem{Remark}[Theorem]{Remark}

\title{Viscosity solutions for a polymer crystal growth model}

\begin{document}

\author[Cardaliaguet, Ley and Monteillet]{Pierre Cardaliaguet 
 \and {Olivier Ley} \and Aur\'elien Monteillet}

\address{
(P. Cardaliaguet, A. Monteillet) Laboratoire de Math\'ematiques \\ CNRS UMR 6205\\  Universit\'e de Brest \\ 6 Av. Le Gorgeu
BP 809, F-29285 Brest, France \\
 {\tt \{pierre.car\-da\-liaguet, aurelien.monteillet\}@univ-brest.fr}} 
\address{ (O. Ley) IRMAR, INSA de Rennes, F-35043 Rennes, France
 \\ {\tt olivier.ley@insa-rennes.fr}}

\begin{abstract}
We prove existence of a solution for a polymer crystal growth model describing the movement of a front $(\Gamma(t))$ evolving with a nonlocal velocity.
In this model the nonlocal velocity is linked to the solution of a heat equation with source $\delta_\Gamma$. The proof relies on new regularity results for 
the eikonal equation, in which the velocity is positive but merely measurable in time and with H\"{o}lder bounds in space. From this result, we deduce \textit{a priori} regularity for the front. On the other hand, under this regularity assumption, we prove bounds and regularity estimates for the solution of the heat equation.
\end{abstract} 

\thanks{This work was partially supported by the ANR (Agence Nationale de la Recherche) through MICA project (ANR-06-BLAN-0082)}
 
 

\maketitle  

%\tableofcontents

\section{Introduction}

The paper is devoted to the analysis of following system of equations: 
\be\label{intro:FPP}
\left\{\begin{array}{lll}
i) & u_t(x,t)=\bar g(v(x,t))|Du(x,t)| \qquad {\rm in }\; \R^N\times (0,+\infty)\\
ii) & v_t(x,t)-\Delta v(x,t)+\kappa\bar g(v(x,t)){\mathcal H}^{N-1}\lfloor\{u(\cdot,t)=0\}=0\\
& \hspace{7cm} {\rm in }\; \R^N\times (0,+\infty)\\
iii) & v(x,0)=v_0(x), \; u(x,0)=u_0(x) \qquad {\rm in }\; \R^N\;.
\end{array}\right.
\ee
Following \cite{ burger01, burger04, bcs02, fv01}, 
the $3$-dimensional version of
this system modelizes the growth of the surface $\Gamma(t)$ of a polymer crystal in a nonhomogeneous temperature field $v(x,t)$. 
In this model one describes the evolving surface $\Gamma(t)$ of the crystal
as the $0$-level-set of an auxiliary function $u$: 
$$\{x\in \R^N\;;\; u(x,t)=0\}=\Gamma(t)\;.
$$ 
(This is the level-set approach, see \cite{giga06} and references therein). 
It has experimentaly been observed that the normal velocity $V_n$ of the crystal is a known, positive function of
the temperature: $V_n=\bar g(v(x,t))$, where $\bar g$ is a bell-shaped function depending on the specific polymer (\cite{Ed97}). 
Expressing the normal velocity $V_n$ in terms of the function $u$ gives the
eikonal equation \eqref{intro:FPP}-$i),$ which holds at least on the set $\{u(\cdot,t)=0\}$.
As for the temperature field $v$ 
it has to follow a heat equation with a (negative) heat source proportional to $V_n{\mathcal H}^{N-1}\lfloor\Gamma(t)$. Whence 
\eqref{intro:FPP}-$ii)$.

Similar systems, coupling eikonal and diffusion equations, 
appear in many applications: shape optimization, image segmentation, etc. (see for instance \cite{Os03, Se99}
and the references therein). However the mathematical analysis of such couplings
is delicate and few existence or uniqueness results are available in the literature. Most of them are concerned with classical
solutions on a short time interval. For instance short time existence and uniqueness of smooth solutions are obtained 
for system \eqref{intro:FPP} in \cite{fv01}.

The point is that, in general, one cannot expect such a system to have classical solutions when the time becomes large: indeed the front $\Gamma(t)$
usually develops singularities in finite time. For this reason a good description of
this front is obtained by its representation as the $0$-level-set of the solution  
of an eikonal equation, which has to 
be understood in the sense of viscosity solutions. However this approach (which is satisfactory from a numerical view point) raises 
severe mathematical difficulties. Such issues have been overcome
in only a very few number of situations: for a dislocation dynamics model, introduced in \cite{ahlm06} and analyzed 
in \cite{acm05, bl06, bclm08}, or for a system arising in the study of the 
asymptotics of a Fitzhugh-Nagumo model \cite{bclm09a, ggi92, ss96}. 
In this later framework, the associated heat equation is of the form
\be\label{Fitz}
v_t(x,t)-\Delta v(x,t)-\bar g(v(x,t)){\bf 1}_{\{u(\cdot,t)\geq 0\}}=0\, ,\\
\ee
where ${\bf 1}_E$ is the indicator function of a set $E$. 
In \cite{bclm09a, ggi92, ss96} existence of generalized solutions for this Fitzhugh-Nagumo system is proved, while
\cite{bclm09} contains some uniqueness results. However, system \eqref{intro:FPP} turns out to be much more challenging than the coupling in the Fitzhuch-Nagumo system. Indeed the surface term
${\mathcal H}^{N-1}\lfloor\{u(\cdot,t)=0\}$ in \eqref{intro:FPP}-$ii)$ is more 
singular than the volume one ${\bf 1}_{\{u(\cdot,t)\geq 0\}}$
in \eqref{Fitz}. For this reason, up to now, only the long time existence in space dimension 
$N=2$ is known \cite{sb07, sb06}.
 
The aim of our paper is to obtain a similar existence result for the physical dimension $N=3$ (and in fact in any dimension). 
In order to state precisely our main result, let us introduce the definition of a solution to \eqref{intro:FPP}. 

\begin{Definition} A solution $(u,v)$ of \eqref{intro:FPP} 
on the time interval $[0,T]$ is a map
$(u,v):\R^N\times [0,T]\to \R^2$ which is  bounded, uniformly continuous, such that
$u$ satisfies the equation 
$$
u_t(x,t)=\bar g(v(x,t))|Du(x,t)|\; {\rm in }\; \R^N\times (0,T),\qquad  u(x,0)=u_0(x)\; {\rm in }\; \R^N
$$
in the viscosity sense, with
$$
\int_0^T {\mathcal H}^{N-1}(\{u(\cdot,t)=0\}) <+\infty\;,
$$
and such that $v(\cdot,0)=v_0$ and $v$ satisfies in the sense of distributions
$$
v_t(x,t)-\Delta v(x,t)+\kappa\bar g(v(x,t)){\mathcal H}^{N-1}\lfloor\{u(\cdot,t)=0\}=0\qquad {\rm in }\; \R^N\times (0,T)\;.
$$
\end{Definition}

We introduce the following set of assumptions, denoted by
{\bf (A)} in the rest of the paper. 
\begin{enumerate}

\vspace{.2cm}

\item[{\bf (A1)}] $\kappa$ is a fixed real number ($\kappa$ is positive in 
the case of a negative heat source and negative otherwise), 
$\bar g:\R^N\to \R$ is Lipschitz continuous, bounded, and there 
exist $A,B>0$ such that 
$$
A\leq \bar g(z)\leq B\qquad {\rm for \ all} \  z\in \R\;.
$$

\vspace{.2cm}

\item[{\bf (A2)}] $v_0:\R^N\to\R$ is Lipschitz continuous and bounded.

\vspace{.2cm}

\item[{\bf (A3)}] $u_0:\R^N\to\R$ is Lipschitz continuous and satisfies $\{ u_0=0 \}=\partial \{u_0>0\}.$ Moreover, we assume that $\{ u_0\geq 0 \}$ is compact and 
has the interior ball property of radius $r_0>0$, that is,
\begin{eqnarray}\label{ibp11}
{\rm For \ all} \  x\in K_0, \ {\rm there \ exists} \  y\in K_0, 
\ {\rm with} \ x\in \overline{B}(y,r_0)\subset K_0\;,
\end{eqnarray}
where $\overline{B}(y,r_0)$ is the closed ball of radius $r_0$ centered at $y$.
\end{enumerate}

Our main result states that, under the above assumptions, system \eqref{intro:FPP} has  a solution. More precisely: 

%%%%%%%%%%%%%
\begin{Theorem}\label{main} Under Assumption {\bf (A)},
for any $T>0$, there exists at least one solution to 
System \eqref{intro:FPP}.
This solution is bounded on $\R^N\times [0,T]$ and satisfies,
for all $x,y\in\R^N,$ $0\leq s, t\leq T,$
$$
|v(x,t)-v(y,t)|\leq C |x-y|(1+|\log |x-y| |),
$$
and
$$
|v(x,t)-v(x,s)|\leq C |t-s|^{\frac12}(1+|\log |t-s| |).
$$
for some constant $C$ which only depends on the data appearing
in Assumption {\bf (A)} and $T.$
\end{Theorem}
%%%%%%%%%%%%%%

Note that  uniqueness of the solution is an open problem (even in dimension~2). 

Let us now briefly describe the method of proof. 
The main difficulty in System \eqref{intro:FPP} is the singular surface
term in the heat equation: 
to deal with this term, one has to obtain fine regularity estimates
for the level-sets of $u.$
Such estimates, which cannot
be derived from the usual regularity results on the eikonal equation, 
have been investigated through several works. 
When the velocity $x\mapsto \bar g(v(x,t))$ is positive of class
${\mathcal C}^{1,1}$, the front enjoys the interior ball property \eqref{ibp11}
\cite{cf06} (see also \cite{acm05, bclm08}); it
has an interior cone property when the velocity is positive and Lipschitz continuous \cite{bclm09}. Unfortunately, for System \eqref{intro:FPP}, 
the interior cone property 
is not sufficient for guarantying the stability of the surface term ${\mathcal H}^{N-1}\lfloor\{u(\cdot,t)=0\}$.
Moreover  we were only able to prove that the map $x\mapsto v(x,t)$ has a modulus of continuity of the form $\omega(\rho)=\rho(1+|\log(\rho)|)$
(even when the front is smooth this map is at most Lipschitz continuous \cite{fv01}). 
Our main and new estimate on the eikonal equation is an interior paraboloid property for the level-sets of $u.$ We call paraboloid
a solid deformation of the set 
$$
\left\{x=(x',x_N)\in \R^{N-1}\times \R
\;;\; x_N\geq c|x'|^{1+\gamma}\right\}, \quad c>0, \ 
\gamma\in (0,1).
$$
This property is obtained under the (weak)
assumption that the velocity $x\mapsto \bar g(v(x,t))$ is of class ${\mathcal C}^{0,\alpha}$. For this, we use a representation 
formula for the solutions of \eqref{intro:FPP}-$i)$ in terms of optimal control as well as sharp regularity properties of  
optimal solutions for this control problem. As a direct consequence of the  interior paraboloid property one obtains
that the front has an interior cone property. These interior paraboloid and cone properties are the two key ingredients
which allow us to obtain
{\it a priori} estimates on the heat flow: indeed, because of the cone property,  the front $\Gamma(t)$ can be covered by a finite
(and controlled) number of Lipschitz graphs. 
The stability result on the surface term ${\mathcal H}^{N-1}\lfloor\{u(\cdot,t)=0\}$ (see 
Lemma \ref{continuite}) is a consequence of the interior paraboloid property. 
Let us finally point out that, although the cone and paraboloid properties do not appear in
 \cite{sb07, sb06}, we use several arguments from these papers: in particular the regularity
of the optimal solutions of some control problem is borrowed from \cite{sb07, sb06} and some of our estimates on the heat flow are related
with those of \cite{sb07, sb06}.

The  paper is organized as follows: Section 2 is dedicated to estimates on the eikonal equation, while the {\it a priori} estimates
for the heat flow are the object of Section 3. We prove the main result in Section 4. \\

\noindent {\bf Notations: } For any integer $k\geq 1$ we denote by $B_k (x,r)$
(resp. $\overline{B}_k(x,r)$) the open (resp. closed) ball of radius
$r>0$ and of center $x$ in $\R^k$. For $k=N$ (the ambiant space), we simply 
abbreviate to $B(x,r)$. We also denote by $\Sc^{N-1}$ the unit sphere of $\R^N$.

%%%%%%%%%%%%%%%%%%%%%%%%%%%%%%%%%%%%%%%%%%%%%%%%%%%%%%%%%%%%%%%%%%%%%%%%%%%%%%%%%%%%%%%%%%%%%%%%%%%%%%%%%%%
\section{Representation formula and a priori estimates for the eikonal equation}

Throughout this section, we investigate the eikonal equation
\begin{equation}\label{Eqeikonal}
\left\{\begin{array}{l}
u_t=c(x,t)|Du|\qquad {\rm in }\; \R^N\times (0,T),\\
u(x,0)=u_0(x) \qquad {\rm in }\;\R^N.
\end{array}\right.
\end{equation}
We assume that the velocity $c$ is Borel measurable on $\R^N \times [0,T]$ and satisfies 
\be \label{BornesVitesse}
A\leq c(x,t)\leq B\qquad {\rm for \ all} \ (x,t)\in \R^N \times [0,T]
\ee
for some $A,B>0$. We also assume that there exist $\alpha\in (0,1)$, $\omega\in L^p(0,T)$ with $p\in (1,+\infty]$ and $C>0$ such that for all $(x,y,t)\in \R^N\times \R^N \times [0,T]$,
\be \label{estiLn}
|c(x,t)-c(y,t)|\le C |y-x| \, (1+|\log |x-y||)\; ,
\ee
and 
\be \label{Holder}
|c(x,t)-c(y,t)|\le \omega(t) |y-x|^\alpha\, .
\ee
Finally, the initial datum $u_0$ is Lipschitz continuous on $\R^N$.
Our aim is to prove existence and uniqueness for the solution 
of \eqref{Eqeikonal} under assumptions
\eqref{BornesVitesse} and \eqref{estiLn}, and give some estimates depending only on 
assumption \eqref{Holder}. Note that the first two parts are quite classical: they are given here for sake of completeness
and also because we are working in a framework (assumption \eqref{estiLn}) which slightly differs from the standard one. 
In constrast, the regularity results on the optimal solutions for the controlled system associated with equation \eqref{Eqeikonal}
and its consequence on the level-sets of the solution of \eqref{Eqeikonal} are new. Their proofs borrow some ideas of 
\cite{sb06,sb07}, as for instance Lemma  \ref{evalby}.  

%%%%%%%%%%%%%%%%%%%%%%%%%%%%%
\subsection{Existence, uniqueness, stability and representation formula}

Let us recall some known results for Equation \eqref{Eqeikonal}. 
 The notion of $L^1$-viscosity solution provides a framework for 
equations such as \eqref{Eqeikonal} where the dependance on the time variable 
is merely measurable. We refer  to \cite[Appendix]{bclm08} for the definition 
and properties of $L^1$-viscosity solutions that we need here, and to
\cite{ishii85, nunziante90, nunziante92, bourgoing04a,bourgoing04b}
for a complete overview of the theory.

Let us introduce the following controlled system: for any $b\in L^\infty([0,T),\R^N)$,
%%%%%%%%%
\be\label{controlsyst}
x'(s)=c(x(s),s)\, b(s) \qquad |b(s)|\le 1, \quad \text{for a.e.} \ s\geq 0.
\ee
%%%%%%%%%

We start by recalling that, for a given initial data and a given control, equation \eqref{controlsyst} has a unique solution
(this is Osgood's Theorem, see \cite{cl55} for instance):
%%%%%%%%
\begin{Lemma} \label{gronwall}
Assume that the function $c:\R^N \times [0,T]\to \R$ is Borel measurable, bounded 
and satisfies \eqref{estiLn}.
For any fixed $b\in L^\infty([0,T),\R^N)$, with $|b(s)|\leq 1$ a.e., Equation 
$$
 \left\lbrace
\begin{aligned}
&x'(s)=c(x(s),s)\, b(s)\quad \text{for a.e.} \ s\in [0,T], \\
&x(0)=x_0 
\end{aligned}
\right.
$$
has a unique absolutely continuous solution on $[0,T]$. Moreover, if $x$ and $y$ are two solutions
of \eqref{controlsyst}, associated to the same control $b\in L^\infty([0,T),\R^N)$, then 
\be \label{form-unicite-edo}
|x(t)-y(t)|\leq \tilde \omega(|x(0)-y(0)|)
\ee
for some modulus $\tilde \omega$ which only depends on the constant $C$ in Assumption 
\eqref{estiLn}. 
\end{Lemma}
\begin{Proposition}[Existence, uniqueness and stability for \eqref{Eqeikonal}] \label{existence}

~~

Assume that the velocity $c:\R^N \times [0,T]\to \R$ is Borel measurable and satisfies \eqref{BornesVitesse} and \eqref{estiLn}. Let $u_0:\R^N\to \R$ be a Lipschitz continuous function. Then:
\begin{itemize}

\item[{\it (i)}] (Existence and uniqueness) Equation 
\eqref{Eqeikonal} has a unique $L^1$-viscosity solution satisfying
\be \label{module}
u_0(x) \leq u(x,t) \leq u_0(x) + B\| Du_0 \|_\infty t\, ,
\ee
for any $(x,t)\in \R^N \times [0,T]$. 

\vspace{.2cm}

\item[{\it (ii)}] (Properties and representation formula)
This solution is nondecreasing in time, 
uniformly continuous on $\R^N \times [0,T]$ and given by the formula
\be \label{representation-u}
u(x,t)=\sup\{ u_0(y) ; \; \exists \, \bar{x} \;\mbox{\rm solution 
of \eqref{controlsyst} with $\bar{x}(0)=y$ and $\bar{x}(t)=x$}\}\, .
\ee
In particular, 
\begin{align} \label{representation}
K(t):&=\{x\in \R^N;\;u(x,t)\geq 0\}\\
&=\left\{x\in \R^N;\; \; \exists \, \bar{x}\;\mbox{\rm solution 
of \eqref{controlsyst} with $\bar{x}(0)\in K(0)$ and $\bar{x}(t)=x$}\right\}. \nonumber
\end{align}

\item[{\it (iii)}] (Stability) If $(c_n)$ is a sequence of measurable functions satisfying 
\eqref{BornesVitesse} and \eqref{estiLn} with the same constants $A,B,C>0$ and such that $(c_n)$ converges a.e. 
to some $c:\R^N\times [0,T]\to\R$, then the sequence of solutions $(u_n)$ of 
\eqref{Eqeikonal} associated to the velocities $(c_n)$ converges locally 
uniformly to the solution $u$ associated to $c$. 

\end{itemize}
\end{Proposition}
%%%%%

{\bf Proof:} The existence of a solution $u$ which satisfies \eqref{module} is a consequence 
of the general theory (see \cite[Propositions 2.1 and 2.2]{nunziante92}). To prove that this 
solution is unique and given by $\eqref{representation-u}$, we proceed by approximation: 
let $(\rho_n)_{n\geq 1}$ be a mollifier on $\R^N$ such that 
${\rm supp}(\rho_n) \subset \overline{B}(0,1/n)$, $\rho_n \geq 0$ and $\| \rho_n \|_{1}=1$. 
Let $(\tilde{c}_n)_{n\geq 1}$ be the sequence of approximate velocities defined by
$$
\tilde{c}_n(x,t)=\int_{\R^N} c(x-y,t) \, \rho_n(y)\, dy.
$$
Then $\tilde{c}_n$ is Borel measurable on $\R^N \times [0,T]$, 
Lipschitz continuous in space (with a $n$-dependant constant), 
satisfies \eqref{BornesVitesse} and \eqref{estiLn}, 
and $(\tilde{c}_n)$ converges to $c$ as $n\to +\infty$. More precisely, 
using \eqref{estiLn}, we have for any $(x,t)\in \R^N \times [0,T],$
$$
| \tilde{c}_n(x,t)-c(x,t) | \leq \int_{\overline{B}(0,1/n)} |c(x-y,t)-c(x,t)|\, 
\rho_n(y) \, dy \leq C\, \frac{1}{n}(1+\log n) .
$$
Let $$c_n^-(x,t)=\tilde{c}_n (x,t)-\frac{C}{n}(1+\log n) \quad
\text{and} \quad c_n^+(x,t)=\tilde{c}_n (x,t)+\frac{C}{n}(1+\log n),$$ so that $c_n^-\leq c\leq c_n^+$ 
and $c_n^{\pm}$ satisfies \eqref{BornesVitesse} with $A/2$ and $2B$ for $n$ large enough. 
By the comparison principle for \eqref{Eqeikonal} with a velocity which is Lipschitz 
continuous in space (see \cite[Theorem 3.1]{nunziante92}), we obtain that 
$u_n^-\leq u\leq u_n^+$, where $u_n^-$ (resp. $u_n^+$) is the solution of 
\eqref{Eqeikonal} associated to the velocity $c_n^-$ (resp. $c_n^+$). 
Moreover \eqref{module} (with $2B$) and  \eqref{representation-u} hold
for both $u_n^-$ and $u_n^+.$  
To conclude, it only remains to prove that, if a sequence of velocities $(c_n)$ 
satisfies \eqref{BornesVitesse} and \eqref{estiLn}, and converges almost everywhere 
to $c$ as $n\to +\infty$, then the representation formulae for the corresponding 
solutions $u_n$ converge to the representation formula for $u$.

First of all, fix $(x,t)\in \R^N\times [0,T]$ and let $(y_n)$ be a sequence of points 
in $\R^N$ such that $u_0(y_n)\to z\in \R$ as $n\to +\infty$ and for any $n$, 
there exists an absolutely continuous function $\bar{x}_n:[0,t]\to \R^N$ 
such that $\bar{x}_n(0)=y_n$, $\bar{x}_n(t)=x$ and 
$|\bar{x}_n'(s)| \leq c_n(\bar{x}_n(s),s)$ on $[0,t]$. Since $|c_n| \leq B$ for any $n$, 
up to an extraction, $(\bar{x}_n)$ converges uniformly to some $\bar{x}:[0,t]\to \R^N$. 
As a consequence, $\bar{x}(t)=x$, $u_0(\bar{x}(0))=z$ and, using the a.e. convergence 
of $(c_n)$ to $c$ as well as \eqref{BornesVitesse} and \eqref{estiLn}, we obtain 
$|\bar{x}'(s)|\leq c(\bar{x}(s),s)$ on $[0,t]$. This proves that
$$
\limsup  u_n(x,t) \leq \sup\{ u_0(y) ; \; \exists \, \bar{x}\;\mbox{\rm solution of 
\eqref{controlsyst} with $\bar{x}(0)=y$ and $\bar{x}(t)=x$}\}\, .
$$

Conversely, let $y\in \R^N$ such that there exists a solution $\bar{x}$
of \eqref{controlsyst} with $\bar{x}(0)=y$ and $\bar{x}(t)=x$. 
Let $b$ be the control associated by $\bar x$ and 
$\bar{x}_n$ be the solution of $\bar{x}_n'(s)=c_n(\bar{x}_n(s),s)b(s)$ with $\bar{x}_n(t)=x.$ 
Then we must have $u_n(x,t)\geq u_0(\bar{x}_n(0))$ for any $n$. By the same argument as above, 
$(\bar{x}_n)$ must converge uniformly to a solution of $x'(s)=c(x(s),s)b(s)$, and by uniqueness 
of such a solution (Lemma \ref{gronwall}), the limit $(x_n)$ must be $\bar{x}$. Therefore
$$
u_0(y) =\lim u_0(\bar{x}_n(0)) \leq \liminf u_n(x,t), 
$$
and
$$
\sup\{ u_0(y) ; \; \exists \, \bar{x}\;\mbox{\rm solution of \eqref{controlsyst} 
with $\bar{x}(0)=y$ and $\bar{x}(t)=x$}\} \leq \liminf  u_n(x,t).
$$
This concludes the proof of the representation formula \eqref{representation-u}
for the unique solution 
of \eqref{Eqeikonal}. This representation formula implies that $u$ is nondecreasing in time. We also point out that the proof of uniqueness can be easily adapted
to prove that, in fact, comparison holds for  \eqref{Eqeikonal}.

\vspace{.2cm}

To prove the stability property $(iii)$, let $(c_n)$ be a sequence of functions satisfying \eqref{BornesVitesse} and \eqref{estiLn} with the same constants $A,B$ and $C$, and such that $(c_n)$ converges a.e. to some $c:\R^N\times [0,T]\to\R$, and let $(u_n)$ be the sequence of solutions of \eqref{Eqeikonal} associated to the velocities $(c_n)$. Using the same arguments as above and the representation formula \eqref{representation-u}, we can actually prove that the half-relaxed limits
$$
{\liminf}_* u_n:(x,t)\mapsto \underset{n\to +\infty}{\liminf}\, \{ u_n(x_n,t_n);\; x_n \to x,\; t_n\to t \}
$$
and
$$
{\limsup}^* u_n:(x,t)\mapsto \underset{n\to +\infty}{\limsup}\, \{ u_n(x_n,t_n);\; x_n \to x,\; t_n\to t \}
$$
coincide and are equal to the solution $u$ of \eqref{Eqeikonal} associated to $c$. This is known to imply the locally uniform convergence of $(u_n)$ to $u$, and proves the stability property.

\vspace{.2cm}

Finally, let us prove the uniform continuity of the solution $u$ of $\eqref{Eqeikonal},$ starting with 
the regularity in space: fix $(x,y,t) \in \R^N \times \R^N \times [0,T]$, and let $\bar{x}$ be 
a solution of \eqref{controlsyst} with control $\bar b$, $\bar{x}(t)=x$ and $u(x,t)=u_0(\bar{x}(0))$
(notice that the supremum is achieved in \eqref{representation-u}). 
Let $\bar{y}$ be the solution of \eqref{controlsyst} associated to the same control $\bar b$ 
and satisfying $\bar{y}(t)=y.$ Applying \eqref{form-unicite-edo} for System \eqref{controlsyst} with reverse
time, we have
$$
|\bar{x}(0)-\bar{y}(0)|
\leq  \tilde \omega (|\bar{x}(t)-\bar{y}(t)|)
 \quad {\rm for \ all} \  t\in [0,T].
$$
Using that $\bar{y}$ is a solution of \eqref{controlsyst} 
and $u_0(\bar{y}(0))\leq u(y,t)$ thanks to \eqref{representation-u}, we obtain
\begin{eqnarray*}
u(x,t)=u_0(\bar{x}(0)) 
&\leq &
u_0(\bar{y}(0)) + \|Du_0\|_\infty \, |\bar{x}(0)-\bar{y}(0)|\\
& \leq &
u(y,t)+ \bar\omega (|x-y|),
\end{eqnarray*} 
where $\bar\omega=\|Du_0\|_\infty \tilde \omega$ is still a modulus of continuity. 
Exchanging the roles of $x$ and $y$, we obtain the uniform continuity 
of $u$ in space. 

Now let us fix $t\in [0,T]$. The map 
$(x,s)\mapsto u(x,t+s)$ is a sub-solution 
of $\bar{u}_t=B\, |D\bar{u}|$ in $\R^N \times [0,T-t]$ with uniformly continuous 
initial datum $u(\cdot,t).$ By the Lax formula, for any $0\leq s\leq T-t$,
$$
u(x,t)\leq u(x,t+s) \leq \sup \{ u(y,t);\; |x-y|\leq Bs \}\, . 
$$
Using the uniform continuity of $u(\cdot,t)$ in space, we deduce that for any $0\leq s\leq T-t$,
$$
u(x,t)\leq u(x,t+s) \leq u(x,t)+\bar\omega (Bs)\, .
$$
This proves the uniform continuity of $u$ in time.

\QED

%%%%%%%%%%%%%%%%%%%%%%%%%%%%%%
%%%%%%%%%%%%%%%%%%%%%%%%%%%%%%
\subsection{Properties of the minimal time function}

Let us now introduce the function
$$
z:x\mapsto \min \{t\in [0,T];\;u(x,t)\geq 0  \}\, ,
$$
which by definition is well-defined on $K(T)=\cup_{t\in [0,T]}\, K(t)$ 
(see \eqref{representation} for the definition of $K(t)$)
and is such that $K(t)=\{x\in \R^N;\; z(x)\leq t\}.$ 

We say that a solution $\bar{x}$ of \eqref{controlsyst} on $[0,t]$ is extremal if 
$$
\bar{x}(0)\in K(0) \quad {\rm and} \quad
z(\bar{x}(t))=t.
$$

%%%%%%%%%%
\begin{Lemma} \label{extremality}
Assume that the velocity $c:\R^N \times [0,T]\to \R$ is Borel measurable and satisfies \eqref{BornesVitesse} and \eqref{estiLn}. 

\begin{enumerate} 
\item Let $\bar{x}$ be an extremal solution on $[0,t]$. Then: 

\vspace{.2cm}

$(i)$ For any $s\in [0,t],$ $z(\bar{x}(s))=s.$

\vspace{.2cm}

$(ii)$ For almost every $s\in [0,t]$, $|\bar{x}'(s)|=c(\bar{x}(s),s)$.\\

\item If $\{x\in\R^N;\; u_0(x)=0\}=\partial \{x\in\R^N;\; u_0(x)>0\}$, then for any $t\in (0,T],$ $$\{ x\in \R^N;\; u(x,t)=0 \}=\{ x\in \R^N;\; z(x)=t \}\, .$$
\end{enumerate}
\end{Lemma}
%%%%%%%%%%%%

{\bf Proof:} $(1)$ $(i)$ By definition of $\bar{x}$ and $z$, we have for any $s\in [0,t]$, 
$z(\bar{x}(s))\leq s$. To prove the converse inequality, we argue by contradiction: 
let $s_0\in [0,t)$ be such that $\theta:=z(\bar{x}(s_0))<s_0.$ Let us first prove that 
for $\delta >0$ small enough,
$$
\overline{B}(\bar{x}(s_0),A(s_0-\theta-\delta))\subset \{y\in \R^N;\; z(y)\leq s_0-\delta \}\, .
$$
Let $y$ be such that $|y-\bar{x}(s_0)|<A(s_0-\theta-\delta),$ 
and let $x_\theta$ be a solution of \eqref{controlsyst} on $[0,\theta]$ such that 
$x_\theta(0)\in K(0)$ 
and $x_\theta(\theta)=\bar{x}(s_0)$. We extend $x_\theta$ to $[0,s_0-\delta]$ by setting
$$
x_\theta(s)=\bar{x}(s_0)+\frac{y-\bar{x}(s_0)}{s_0-\theta-\delta}\, 
(s-\theta) \quad {\rm for \ all} \  s \in [\theta,s_0-\delta]\, .
$$
The bound $c \geq A$ shows that $x_\theta$ is a solution of \eqref{controlsyst} 
on $[0,s_0-\delta]$ with $x_\theta(0)\in K(0)$ and $x_\theta(s_0-\delta)=y,$ 
which means that $z(y)=z(x_\theta(s_0-\delta))\leq s_0-\delta.$

Now, for any $\delta >0$ small enough, let us solve
$$
\left\lbrace
\begin{aligned}
&x_{\delta}'(s)=c(x_\delta(s),s)\, b(s) \quad \text{on} \; [s_0-\delta,t-\delta]\, ,\\
&x_\delta(t-\delta)=\bar{x}(t)\, .
\end{aligned} \right.
$$
where $b$ is the control associated to $\bar{x}$. 
Applying \eqref{form-unicite-edo} for System \eqref{controlsyst} with reverse
time, we have
\begin{eqnarray*}
|x_\delta(s_0-\delta)-\bar{x}(s_0-\delta)|
&\leq& \tilde \omega(|x_\delta(t-\delta)-\bar{x}(t-\delta)|)\\
&=&\tilde 
\omega(|\bar{x}(t)-\bar{x}(t-\delta)|)\\
&\leq &\tilde 
\omega(B\delta)
\end{eqnarray*}
because $|\bar{x}'|\leq B.$ In particular, for $\delta$ small enough,
$$
|x_\delta(s_0-\delta)-\bar{x}(s_0-\delta)| < \frac12 A(s_0-\theta-\delta)\, ,
$$
while
$$
|\bar{x}(s_0-\delta)-\bar{x}(s_0)| \leq B\delta < \frac12 A(s_0-\theta-\delta)\, .
$$
For such a choice of $\delta$, 
$$
x_\delta(s_0-\delta)\in \overline{B}(\bar{x}(s_0),A(s_0-\theta-\delta))
\subset \{y\in \R^N;\; z(y)\leq s_0-\delta \}.
$$ 
Therefore $z(x_\delta(s_0-\delta))\leq s_0-\delta.$ 
In particular, there exists a solution $\tilde{x}$ of \eqref{controlsyst} 
on $[0,s_0-\delta]$ with $\tilde{x}(0)\in K(0)$ and 
$\tilde{x}(s_0-\delta)=x_\delta(s_0-\delta).$ The reunion of the paths associated 
to $\tilde{x}$ on $[0,s_0-\delta]$ and $x_{\delta}$ on $[s_0-\delta,t-\delta]$ 
gives a solution $x$ of \eqref{controlsyst} on $[0,t-\delta]$ 
with $x(0)\in K(0)$ and $x(t-\delta)=x_\delta(t-\delta)=\bar{x}(t).$ 
In particular, $z(\bar{x}(t))\leq t-\delta <t$, which is absurd.

\vspace{.2cm}

$(1)$ $(ii)$ Now, let us prove that $|\bar{x}'(s)|=c(\bar{x}(s),s)$ 
for almost every $s\in [0,t]$: 
for $s_0\in(0,t)$ and $h>0$ be small enough, 
let $y:[s_0-h,s_0+h]$ be the solution of
$$
\left\{\begin{array}{l}
y'(s)=c(y(s),s)\frac{\bar{x}(s_0+h)-\bar{x}(s_0-h)}{|\bar{x}(s_0+h)-\bar{x}(s_0-h)|}\, ,\\
y(s_0-h)=\bar{x}(s_0-h)\, .
\end{array}\right.
$$
($\bar x$ is injective from (1) $(i)$). 
Note that $y$ remains in the segment $[\bx(s_0-h),\bx(s_0+h)]$ on $[s_0-h,s_0+h]$ 
because $|y'(s)|\leq c(y(s),s)$, which means that $y$ is sub-optimal. 
Moreover $y$ is monotonous on this segment. In particular we have
$$
|\bx(s_0+h)-\bx(s_0-h)| \geq |y(s_0+h)-y(s_0-h)|=\int_{s_0-h}^{s_0+h} c(y(s),s)\, ds.
$$
Using the bound $c\leq B$, we have 
$$
|y(s)-\bx(s)|\leq 4B h \quad {\rm for \ all} \ s\in [s_0-h,s_0+h]\;.
$$
Therefore, thanks to \eqref{estiLn}, we get
$$
\int_{s_0-h}^{s_0+h} c(y(s),s)\, ds \geq 
\int_{s_0-h}^{s_0+h} c(\bx(s),s)\, ds - 8BC h^2 (1+| \log (4B h)|).
$$
If $s_0$ is a Lebesgue point of $s \mapsto c(\bar{x}(s),s)$ such that $\bar{x}$ is differentiable at $s_0$, which is the case of almost every $s_0\in[0,t]$, then we obtain
$$
\begin{aligned}
| \bx'(s_0)| &=\underset{h\to 0}{\lim}\; \frac{|\bx(s_0+h)-\bx(s_0-h)|}{2h}\\ &\geq \underset{h\to 0}{\lim}\; \frac{1}{2h} \int_{s_0-h}^{s_0+h} c(\bx(s),s)\, ds =c(\bx(s_0),s_0)\, .
\end{aligned}
$$

\vspace{.2cm}

$(2)$ Let $(x,t)\in \R^N \times (0,T]$ be such that $z(x)=t;$ by definition of $z$, 
we know that $u(x,t)\geq 0$ and for any $h>0$ enough, $u(x,t-h)<0$. By continuity of $u$, 
we must have $u(x,t)=0.$

Conversely, let $(x,t)\in \R^N \times (0,T]$ be such that $u(x,t)=0.$ We argue by 
contradiction and assume that $\theta=z(x)<t.$ 
Since $u$ is nondecreasing in $t,$ one necessarily has $u(x,\theta)=0.$
Let $\bx$ be a solution of 
\eqref{controlsyst} such that $u(x,\theta)=u_0(\bx(0))=0$ and $\bx(\theta)=x.$ 
By our assumption on $u_0$, there exists $y$ such that $u_0(y)>0$ and 
$$
\tilde \omega (|y-\bar{x}(0)|) < A(t-\theta)
$$ 
(recall that $\tilde \omega$ is defined in \eqref{form-unicite-edo}).
Let $\bar y$ be the solution of \eqref{controlsyst} on $[0,\theta]$ with the control $b$ 
associated to $\bx$, and such that $\bar{y}(0)=y.$ 
Then, from  \eqref{form-unicite-edo}, we have
$$
|\bar{y}(\theta)-x| 
=  |\bar{y}(\theta)-\bar x(\theta)|
\leq \tilde \omega (|\bar{y}(0)-\bar{x}(0)|)<  A(t-\theta).
$$
We extend $\bar{y}$ to $[0,t]$ by setting for any $s\in[\theta,t],$
$$
\bar y(s)=\bar y(\theta)+\frac{x-\bar y(\theta)}{t-\theta}\, (s-\theta)\, .
$$
The bound $c\geq A$ implies that $\bar{y}$ is a solution of \eqref{controlsyst} with $\bar{y}(0)=y$ and $\bar{y}(t)=x.$ By \eqref{representation-u}, we have $u(x,t)\geq u_0(\bar{y}(0))=u_0(y)>0,$ which is absurd. Therefore $z(x)=t$, and this concludes the proof.
\QED

%%%%%%%
\begin{Proposition} \label{zLipschitz}
Under the assumptions of Proposition \ref{existence}, the map $z$ satisfies
$$
\frac1B \leq |Dz|\leq \frac1A
$$
in the viscosity sense and therefore almost everywhere in $\{ x\in \R^N;\; 0<z(x)<T\}.$
\end{Proposition}
%%%%%%%

{\bf Proof:}
The proof of the right-hand side inequality follows along the same lines as the beginning of the proof of \cite[Theorem 5.9]{bclm09}, and shows that $z$ is Lipschitz continuous. For the left-hand side inequality, let $\phi : \{ x\in \R^N;\; 0<z(x)<T\}\to \R$ 
be a function of class $C^1$ such that $z-\phi$ has a local minimum equal to 0 at some $x$. Let $\bar{x}$ be an extremal on $[0,t]$ with $\bar{x}(t)=x$. For any $s\in [0,t]$, $z(\bar{x}(s))=s$ by Lemma \ref{extremality}. Then for any $h>0$ small enough,
$$
z(\bar{x}(t-h))\geq \phi(\bar{x}(t-h))\, ,
$$
whence, by definition of $\phi$,
$$
\phi(\bar{x}(t))-h=z(\bar{x}(t))-h=t-h=z(\bar{x}(t-h))\geq \phi(\bar{x}(t-h))\, .
$$
In particular,
$$
h\leq \phi(\bar{x}(t))-\phi(\bar{x}(t-h))=\int_{t-h}^t \langle D\phi(\bar{x}(s)), \bar{x}'(s)\rangle \, ds \leq B\, \int_{t-h}^t | D\phi(\bar{x}(s))| \, ds
$$
thanks to the the bound $c\leq B$. Dividing this expression by $h$ and letting $h\to 0$, we get $|D\phi(x)|\geq 1/B.$ Since $z$ is Lipschitz continuous, the viscosity inequality $|Dz|\geq 1/B$ also holds almost everywhere.
\QED

%%%%%%%%
\begin{Remark} \label{A consequence}{\rm 
 A consequence of the inequality $|Dz|\geq 1/B$ and Lemma \ref{extremality} (2) is that for any $t\in [0,T]$, the front $\{ x\in \R^N;\;u(x,t)=0 \}$ has measure 0 and coincides with $\partial K(t).$ Indeed, $\{ x\in \R^N;\;u(x,t)=0 \}=\{ x\in \R^N;\;z(x)= t \}$, and Stampacchia's theorem (see for instance \cite{eg92}) states that $Dz=0$ almost everywhere on the set $\{x\in \R^N;\; z(x)=t\}.$ Moreover, the viscosity decrease principle (see \cite{ley01}) shows that
$$
\begin{aligned}
\partial K(t) &=\partial \{ x\in \R^N;\;z(x)\leq t \}= \{ x\in \R^N;\;z(x)= t \}= \{ x\in \R^N;\;u(x,t)=0 \}\, .
\end{aligned}
$$
In particular, a solution $\bar{x}$ of \eqref{controlsyst} is extremal on $[0,t ]$ if 
$x(t )\in \partial K( t )$; in this case, it satisfies $\bar{x}(s)\in \partial K(s)$ for any $s\in [0,t]$
and $|\bar{x}'(s)|=c(\bx(s),s)$ for a.e. $s\in [0,t]$. 
}\end{Remark}
%%%%%%%%

%%%%%%%%%%%%%%%%%%%%%%%%%%%%%%%%%%%%%%%%%%%%%%%%%%%%%%%%%%%%%%%%%
\subsection{Regularity of extremal solutions}

From now on we assume that $c$ satisfies \eqref{BornesVitesse},  \eqref{estiLn} and
\eqref{Holder}. 
Our first result is the following:

%%%%%%%
\begin{Proposition}\label{prop:regux} Under the above assumptions, if $\bx$ is extremal on $[0,\bar t \,]$ for some
$\bar t \,\in(0,T]$ and if $\beta:=\alpha-1/p>0$, then
the map $t\to \bx'(t)/|\bx'(t)|$ is  of class ${\mathcal C}^{\beta/2}(0,\bar t \, )$. Namely
$$
\left|\frac{\bx'(s_2)}{|\bx'(s_2)|}-\frac{\bx'(s_1)}{|\bx'(s_1)|}\right|\leq C\|\omega\|_p^{1/2}\, |s_2-s_1|^{\beta/2} \qquad for \ all \ s_1,s_2\in [0,\bar t \,]\;,
$$
where $C$ only depends on the constants $A,B,\alpha$ and $p$ introduced in \eqref{BornesVitesse}--\eqref{Holder}. 
\end{Proposition}
%%%%%%%

{\bf Proof:} Throughout the proof $C$ denotes a constant which depends on $A,B,\alpha$ and $p$ only. 

By Lemma \ref{extremality} (1)$(ii)$, we have $|\bx'(t)|=c(\bx(t),t)$ a.e. on $[0,\bar t]$. 
We reparametrize the path $\bar x$ with speed 1 as follows.
Let $\theta$ be a solution of 
\be \label{reparam}
\left\{\begin{array}{l}
\theta'(s)=\frac{1}{c(\bx(\theta(s)),\theta(s))} \qquad s\in [0,\theta^{-1}(\bar t \, ) ],\\
\theta(0)=0.
\end{array}\right.
\ee
Let us set $\bar s=\theta^{-1}(\bar t \, )$ and $\by(s)=\bx(\theta(s))$ on $[0,\bar s ]$. Then 
\be \label{vit-y-1}
|\by'(s)|=|\bx'(\theta(s))\theta'(s)|=1  \quad \text{for any} \; s\in [0,\bar s ]\;.
\ee
Let us introduce
\be \label{defcbarre}
\bar c(y,s)=\frac{c(y,\theta(s))}{c(\by(s)),\theta(s))}.
\ee
From our assumptions
\eqref{BornesVitesse}--\eqref{Holder}, we have
\be\label{EstiBarc}
\left|\bar c(y,s)-\bar c(y',s)\right|\leq \frac{\omega(\theta(s))}{A} |y-y'|^\alpha
\quad {\rm for \ all} \ (y,y',s)\in \R^N\times \R^N \times [0,\bar s ]
\ee
and
\be \label{boundbarc}
\frac{A}{B}\leq \bar c(y,s)\leq \frac{B}{A} \quad {\rm for \ all} \ (y,s)\in\R^N \times [0,\bar s ].
\ee

In order to proceed we need the following lemma:
%%%%%%
\begin{Lemma}\label{evalby} For any $0\leq s_1<s_2\leq \bar s$,
$$
|\by(s_2)-\by(s_1)|
\leq s_2-s_1=\int_{s_1}^{s_2}|\by'(s)|ds 
\leq |\by(s_2)-\by(s_1)|+ C(s_2-s_1)^\alpha\int_{s_1}^{s_2}\omega(\theta(s))ds\;.
$$
\end{Lemma}
%%%%%%

{\bf Proof:} First of all, $|\by(s_2)-\by(s_1)|
\leq s_2-s_1=\int_{s_1}^{s_2}|\by'(s)|ds$ because $|\by '|= 1$. Let $y:[s_1,s_2]\to \R^N$ solve
\be \label{systeme-z}
\left\{\begin{array}{l}
y'(s)=\bar c(y(s),s)\frac{\bar y(s_2)-\bar y(s_1)}{|\bar y(s_2)-\bar y(s_1)|}\; , \\
y(s_1)=\bar y(s_1)\;.
\end{array}\right.
\ee
Note that $y$ remains in the segment $[\bar y(s_1),\bar y(s_2)]$ on $[s_1,s_2]$ because $y$
is admissible for \eqref{controlsyst}, and so is sub-optimal. Moreover $y$ is monotonous on the segment. 
From the bounds \eqref{boundbarc} on $\bar c$, we have 
$$
|y(s)-\by(s)|\leq \frac{2B}{A}(s_2-s_1) \quad {\rm for \ all} \  s\in [s_1,s_2]\;.
$$
Since $\bar c(\by(s),s)=1$ and $\bar c$ satisfies \eqref{EstiBarc}, we have
$$
s_2-s_1=\int_{s_1}^{s_2}\bar c(\by(s),s)dt \leq \int_{s_1}^{s_2}\bar c(y(s),s)dt
+ \left(\frac{2B}{A}\right)^\alpha(s_2-s_1)^\alpha
\int_{s_1}^{s_2}\frac{\omega(\theta(s))}{A}ds.
$$
On the other hand, $y$ lives in the segment $[\by(s_1),\by(s_2)]$ and is monotonous on this segment, so that,
from \eqref{systeme-z}, we get
$$
\int_{s_1}^{s_2}\bar c(y(s),s)ds= \int_{s_1}^{s_2}|y'(s)|ds=|y(s_2)-y(s_1)|\leq |\by(s_2)-\by(s_1)|\;.
$$
Putting together the last two estimates proves the Lemma.
\QED

Next we claim the following result:
%%%%
\begin{Lemma}\label{evalby2} For any $0\leq s_1<s_2\leq \bar s$, we have
$$
\begin{array}{l}
\displaystyle{ \left|\by\left(\frac{s_1+s_2}{2}\right)-\frac{\bar y(s_1)+\bar y(s_2)}{2}\right|   }\\
\qquad \qquad \leq 
\displaystyle{ 
C \left\{(s_2-s_1)^\alpha\int_{s_1}^{s_2}\omega(\theta(s))ds+
(s_2-s_1)^{(1+\alpha)/2}\left(\int_{s_1}^{s_2}\omega(\theta(s))ds\right)^{\frac12}\right\}   }.
\end{array}
$$
\end{Lemma}
%%%%

{\bf Proof:} Let us set $s_0=(s_1+s_2)/2$, $a=\by(s_0)-\by(s_1)$, $b=\by(s_2)-\by(s_0)$ and $\tau=s_2-s_1$. Then, from Lemma 
\ref{evalby} we have
$$
|a|+|b|\; \leq \int_{s_1}^{s_0} |\by'(s)|\, ds+\int_{s_0}^{s_2} |\by'(s)|\, ds \; \leq \; \int_{s_1}^{s_2}|\by'(s)|ds \; \leq \;  |a+b|+\ep\; ,
$$
where $\ep:=C\tau^\alpha \int_{s_1}^{s_2}\omega(\theta(s))ds$. 
Taking the square in the above inequality and expanding this expression, we get
$$
2|a||b|\leq 2\lg a,b\rg+2|a+b|\ep+\ep^2\;.
$$
Hence
$$
\left|\frac{a}{|a|}-\frac{b}{|b|}\right|^2 \leq \frac{2|a+b|\ep+\ep^2}{|a||b|}\;.
$$
From \eqref{vit-y-1} and \eqref{boundbarc}, we have 
$$
\frac{A}{B}\frac{\tau}{2} \leq |a|,|b| \leq \frac{\tau}{2}\;.
$$
It follows that
$$
\left|\frac{a}{|a|}-\frac{b}{|b|}\right|^2 \leq 8\left(\frac{B}{A}\right)\frac{\ep}{\tau}
+4\left(\frac{B}{A}\right)^2\frac{\ep^2}{\tau^2}\;.
$$
Let us estimate $||a|-|b||$: from Lemma \ref{evalby} we have
$$
|a|\leq \int_{s_1}^{s_0} |\by'(s)|ds =\frac{\tau}{2}=\int_{s_0}^{s_2} |\by'(s)|ds \leq 
|\by(s_2)-\by(s_0)|+ \ep= |b|+\ep \,.
$$
We obtain the inequality $|b|\leq |a|+\ep$ in the same way, which proves that $||a|-|b||\leq \ep$. Then we write
$$
|a-b|= |a| \left|\frac{a}{|a|}-\frac{b}{|a|}\right|\leq |a| \left|\frac{a}{|a|}-\frac{b}{|b|}\right|+||a|-|b||\;.
$$
Therefore, since $|a|\leq \tau/2$, we have
$$
|a-b|\leq C(\sqrt{\ep \tau}+\ep)\;,
$$
which is the desired result from the definition of $\ep$. 
\QED

We are now ready to complete the proof of Proposition \ref{prop:regux}. 
Since $1/B\leq \theta'\leq 1/A,$ we have
$$
\int_{s_1}^{s_2} \omega(\theta(s))ds= 
\int_{\theta(s_1)}^{\theta(s_2)} \frac{\omega(s)}{\theta'(\theta^{-1}(s))}ds\leq 
B  \int_{\theta(s_1)}^{\theta(s_2)} \omega(s) ds
$$
where, from H\"{o}lder's inequality, 
$$
\int_{\theta(s_1)}^{\theta(s_2)} \omega(s) ds
\; \leq \; |\theta(s_2)-\theta(s_1)|^{1-1/p}\|\omega\|_{p} \;  \leq \;  A^{-1+1/p}|s_2-s_1|^{1-1/p}\|\omega\|_p.
$$
This shows that
\be\label{omegat1t2}
\int_{s_1}^{s_2} \omega(\theta(s))ds \leq C \|\omega\|_p\, |s_2-s_1|^{1-1/p} \;.
\ee
If $\beta=\alpha-1/p>0$, then, combining Lemma \ref{evalby2} with \eqref{omegat1t2}, we get
$$
\left|\by\left(\frac{s_1+s_2}{2}\right)-\frac{\bar y(s_1)+\bar y(s_2)}{2}\right| 
\;  \leq  \; 
C\|\omega\|_p^{1/2} \, |s_2-s_1|^{1+\beta/2}
$$
as soon as $s_2-s_1\leq \|\omega \|_p^{-1/\beta}.$
Theorem 2.1.10 of \cite{cs04} then states that each component of $\by$ is semi-convex and 
semi-concave with a modulus $m$ of the form $m(\rho)=C\|\omega\|_p^{1/2}\rho^{\beta/2}$. Moreover, from Theorem 3.3.7 of \cite{cs04},
we know that $\by$ is ${\mathcal C}^{1,\beta/2}$ with constant $C\|\omega\|_p^{1/2}$. Therefore
$$
|\by'(s_2)-\by'(s_1)|\leq C\|\omega\|_p^{1/2}\, |s_2-s_1|^{\beta/2}
$$
which completes the proof since $\theta^{-1}$ is $B-$Lipschitz continuous and $$\frac{\bx'(t)}{|\bx'(t)|}=\by'(\theta^{-1}(t)).$$
\QED

%%%
\begin{Remark}\label{rem:y} 
\rm We have actually proved that $\by$ is ${\mathcal C}^{1,\beta/2},$
$\beta =\alpha -1/p,$ with constant 
$C\|\omega\|_p^{1/2},$ where $C$ depends only on $A,B,\alpha$ and $p$.
\end{Remark}
%%%

%%%%%%%%%%%%%%%%%%%%%%%%%%%%%%%%%%%%%%%%%%%%%%
\subsection{A priori regularity of the moving front}

We consider a solution  $u$ to \eqref{Eqeikonal} for a velocity $c$ 
which satisfies \eqref{BornesVitesse}, \eqref{estiLn} and \eqref{Holder}. 
 We set, as before, 
$$
K(t)=\{x\in \R^N\,;\, u(x,t)\geq 0\}\quad {\rm for \ all} \  t\in [0,T]\;.
$$
We introduce cone-like sets and interior cone properties as follows.

%%%%%%%
\begin{Definition} \label{def-cones}
Let $x\in \R^N$ and $\nu\in\Sc^{N-1}$ be a unit vector.

\begin{itemize}

\item For any $0<\rho < \theta,$
the {\em cone} of vertex $x,$ axis $\nu$ and parameters $(\rho,\theta)$
is defined by
\begin{eqnarray*} 
\cone_{\nu,x}^{\rho,\theta} 
& := &
\displaystyle \bigcup_{t\in [0,\theta]} B\left(x+t\nu , t\frac{\rho}{\theta} \right)\\
& = & \{x+t\nu+t\frac{\rho}{\theta}\xi\, : \, 
t\in [0,\theta], \, \xi\in\overline{B}(0,1)\}.
\end{eqnarray*}

\item For $C>0$, $\delta\in(0,1)$, we define the {\em paraboloid}
\begin{eqnarray*}
\conereg^{\delta,C}(x, \nu)
&=&
\bigcup_{t\in [0,C^{-1/\delta}]} B\left(x+t\nu , t-Ct^{1+\delta}\right)\\
&=&
\{ x+t\nu+(t-Ct^{1+\delta})\xi\, : \, 
t\in [0,C^{-1/\delta}], \, \xi\in\overline{B}(0,1)
\}.
\end{eqnarray*}

\end{itemize}

We recall from \cite{bclm09} that a compact subset $K$ of $\R^N$
is said to have the 
interior cone property of parameters $(\rho,\theta)$ if, for
any $x\in  \partial K,$ there exists $\nu\in \Sc^{N-1}$ such that the cone
$\cone_{\nu,x}^{\rho,\theta}$ is contained in $K.$

In the same way, we say that $K$ satisfies the interior 
$\conereg^{\delta,C}$-property if for any
$x\in \partial K$, there exists $\nu\in \Sc^{N-1}$ such that
$\conereg^{\delta,C}(x, \nu)$ is contained in $K$.

\end{Definition}
%%%%

\begin{figure}[ht]  
\begin{center}  
\hspace*{-1.5cm}
\epsfig{file=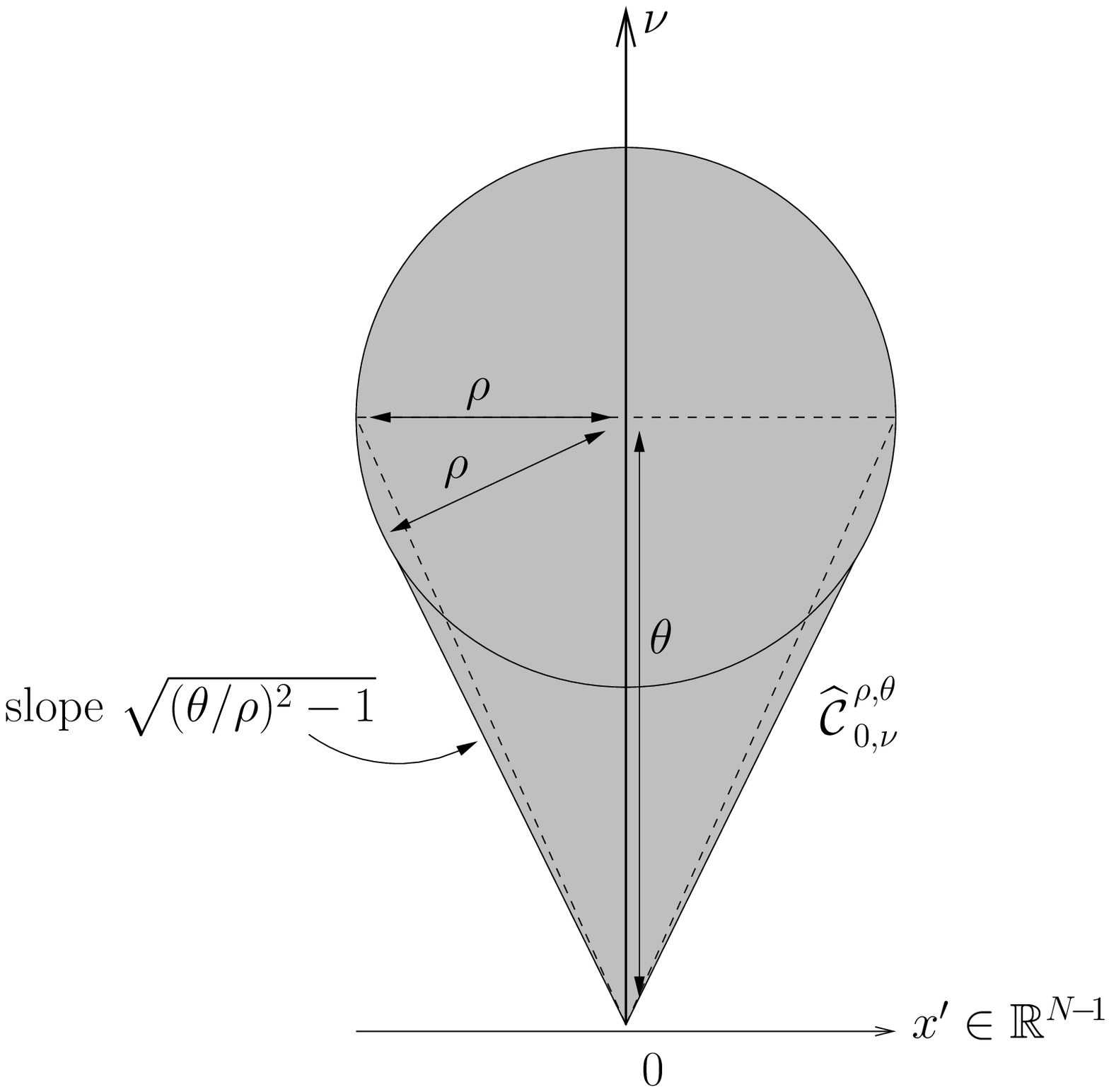, width=7cm}  
\epsfig{file=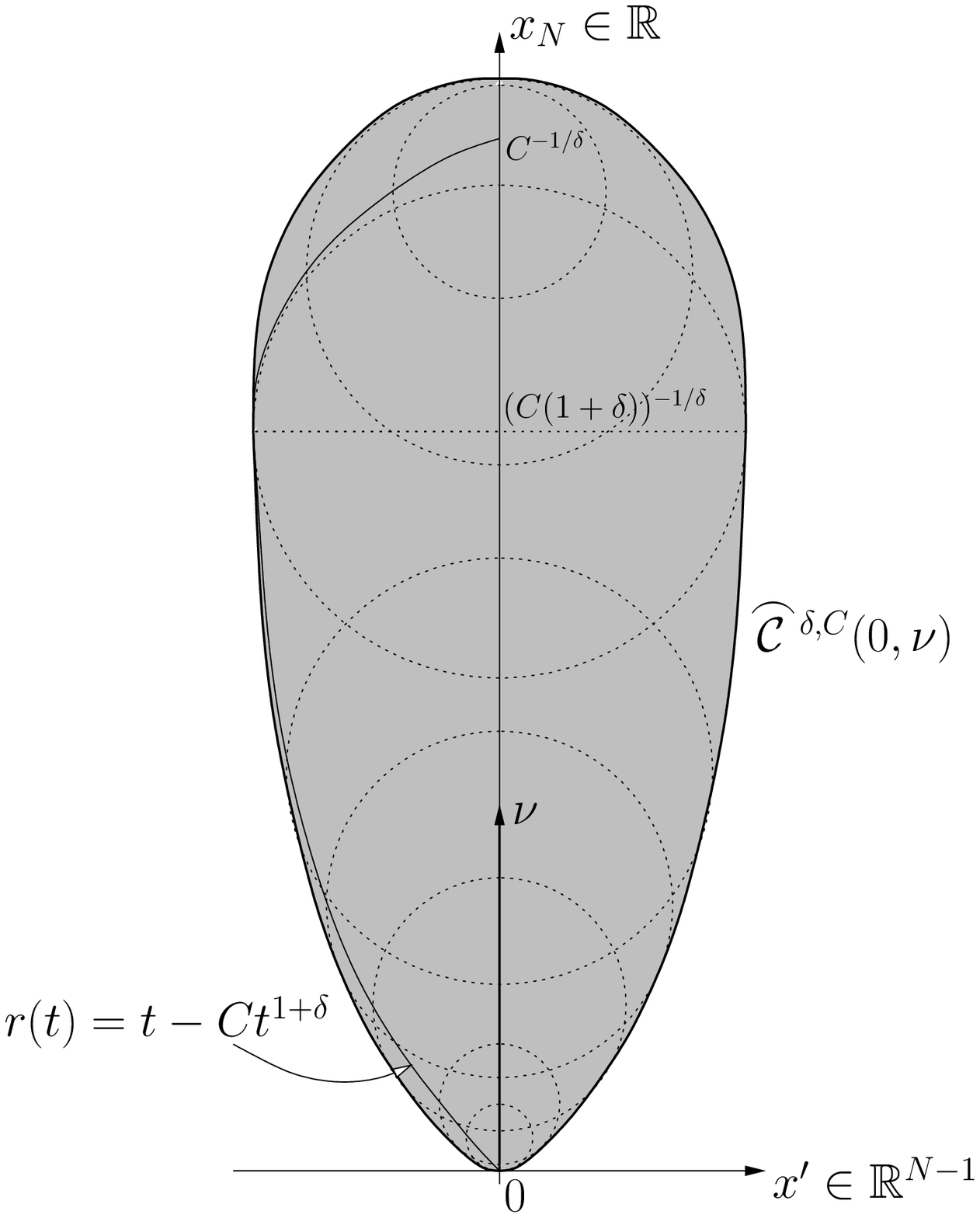, width=6cm}  
\end{center}  
\caption{ \label{dessins-cones}  
{\it Classical cone and paraboloid.
}}  
\end{figure}

The set $\cone_{\nu,x}^{\rho,\theta}$ is a classical cone (see Figure \ref{dessins-cones}).
Since the map $t\to t-Ct^{1+\delta}$ is concave, a tedious but straightforward 
computation shows that the set $\conereg^{\delta,C}(x, \nu)$ is convex.
We shall see below (Lemma \ref{lem:Cgraphe}) that it has a  
${\mathcal C}^{1, \gamma}$ boundary in a neighbourhood
of $x$ for some $\gamma\in (0,1)$ and contains a paraboloid-like subset.
This motivates the name {\em paraboloid}
(see Figure \ref{dessins-cones} 
for an illustration). 
Notice that $\cone_{\nu,x}^{\rho,\theta} \subset \; \conereg^{\delta,C}(x, \nu)$
as soon as $\theta \leq C^{-1/\delta}$ and $\rho \leq \theta - C \theta^{1+\delta}.$

%%%%%%%%%%%%%%%%%%%%%%%%%%%%  
% figure: ensemble cone regulier
%%%%%%%%%%%%%%%%%%%%%%%%%%%%%%%%%%%%%%%%%%%%%%%%%% 

%%%%%%
\begin{Lemma}\label{lem:ccone} Let us still assume that $\beta=\alpha-1/p>0$. 
There exist positive constants $C_0, C_1$ depending only on $A$, $B$, $\alpha$ and $p$,
such that, setting $C(\omega)=C_0 \|\omega\|_p^{1/2},$  
for any extremal solution $\bx$ 
on $[0,\bar t \,]$ with $\bar t \geq C_1 C(\omega)^{-2/\beta}$, 
the set $\conereg^{\beta/2,C(\omega)}(x, \nu)$ is contained in
$K(\bar t \, )$, where
$$
x=\bar x(\bar t \, )\;, \quad \nu=-\frac{\bx'(\bar t \, )}{|\bx'(\bar t \, )|}.
$$
\end{Lemma}
%%%%%%

{\bf Proof:} As in the proof of Proposition~\ref{prop:regux},
we reparametrize $\bar x$ with speed 1 by introducing
$\by(s)=\bx(\theta(s))$ on $[0,\bar s ]$ where $\bar y$ and
$\bar s$ are defined by \eqref{reparam}. 
Notice that $\by'(s)=\bx'(\theta (s) )/|\bx'(\theta (s) )|$ for
a.e. $s\in [0, \bar s ].$ 

Next we define $\bar c$ by \eqref{defcbarre} and,
for $s\in (0,\bar s)$ and $b\in \overline{B}(0,1)$,
we consider the solution
$y:[s,\bar s ]\to \R^N$ to
$$
\left\{\begin{array}{l}
y'(\sigma)=\bar c(y(\sigma),\sigma)b, \qquad \sigma \in [s,\bar s ],\\
y(s)=\by(s).
\end{array}\right.
$$
Arguing as in the proof of Lemma~\ref{evalby}, we obtain
that $y$ is sub-optimal and monotonous on the
segment $[\by(s),y(\bar s)].$
In particular, this whole segment lies in $K(\bar t \, )$.

From the bound \eqref{boundbarc} on $\bar c$, we have
$$
|y(\sigma)-\bar y(\sigma)|\leq \frac{2B}{A}(\bar s-s) 
\quad {\rm for \ all} \ \sigma \in [s,\bar s ].
$$
Hence, by \eqref{EstiBarc},
$$
\bar s-s =\int_{s}^{\bar s} \bar c(\by(\sigma),\sigma)d\sigma \leq \int_{s}^{\bar s} \bar c(y(\sigma),\sigma)d\sigma+ 
\frac{(2B)^\alpha}{A^{1+\alpha}}(\bar s-s)^\alpha \int_s^{\bar s}\omega(\theta(\sigma))d\sigma \; ,
$$
where
$$
\int_s^{\bar s}\omega(\theta(\sigma))d\sigma \leq BA^{-1+1/p}\, \|\omega\|_p\ (\bar s-s)^{1-1/p}\;.
$$
Since $y$ lives in the segment 
$[\by(s),y(\bar s)]\subset [\by(s),\by(s)+\frac{B}{A}(\bar s-s)b]$ 
and is monotonous on this segment, we have
$$
\int_{s}^{\bar s} \bar c(y(\sigma),\sigma)d\sigma= |y(\bar s)-\by(s)|\;.
$$
It follows that
$$
|y(\bar s)-\by(s)|\geq (\bar s-s)(1-\tilde C\|\omega\|_p(\bar s-s)^{\beta}),
\; \text{where} \;
\tilde C =\frac{2^\alpha B^{1+\alpha}}{A^{2+\alpha-1/p}}
\; \text{and} \;
\beta=\alpha-\frac{1}{p}.
$$
Moreover, any point in the segment $[\by(s),y(\bar s)]$ also belongs to $K(\bar t \, )$. 
We have therefore proved that
$$
\by(s)+(\bar s-s)(1-\tilde C\|\omega\|_p(\bar s-s)^{\beta})b\in K(\bar t \, )\;.
$$
This holds true for any $b\in \overline{B}(0,1)$ and any $s$ such that 
$(\bar s-s)\leq \tilde{C}^{-1/\beta}\|\omega\|_p^{-1/\beta}.$
In particular, as soon as $\bar s \geq  \tilde{C}^{-1/\beta}\|\omega\|_p^{-1/\beta},$ we have, setting $t=\bar s - s,$
$$
\bigcup_{t\in [0,  \tilde{C}^{-1/\beta}\|\omega\|_p^{-1/\beta}]} 
B\left(\by(\bar s- t), t(1-\tilde C\|\omega\|_pt^{\beta})\right)   \subset  
K(\bar t \, ) \; ,
$$
where $\beta=\alpha-1/p>0$. From the ${\mathcal C}^{1, \beta/2}$ regularity of $\by$ (see Remark \ref{rem:y}), 
using that $\nu=-\dfrac{\bx'(\bar t \, )}{|\bx'(\bar t \, )|}=-\by'(\bar s),$
we have 
$$
|\by(\bar s-t)-(x+t\nu )|\leq C\|\omega\|_p^{1/2}\int_{0}^{t} s^{\beta/2}ds \leq 
C\|\omega\|_p^{1/2}\, t^{1+\beta/2},
$$
where $C$ only depends on $A,B,\alpha$ and $p$. 
Let us set $$C_0= C+B^{\beta/2}\tilde{C}^{1/2}, \quad C_1=A^{-1}\tilde{C}^{-1/\beta}C_0^{2/\beta}$$
and $$C(\omega)=C_0\|\omega\|_p^{1/2}.$$
Then, going back to the expression of $\bar x,$ we obtain that,
if $\bar t \geq C_1 C(\omega)^{-2/\beta}$,
$$
\conereg^{\beta/2,C(\omega)}(x, \nu)
=\bigcup_{t\in [0, C(\omega)^{-2/\beta}]} B\left(x+ t\nu, t(1-C(\omega)t^{\beta/2})\right)  
\subset K(\bar t \, ).
$$
\QED

The above results have the following consequence:

%%%%%%%%%%%%%
\begin{Corollary}\label{cor:ccone} Let us assume that $K_0$ has the interior ball property of radius $r_0$:
\be\label{IBC}
For \ all \ x\in K_0, \ there \ exists \ y\in K_0, 
\ with \  x\in \overline{B}(y,r_0)\subset K_0\;.
\ee
Then there is a positive constant $C_0$ depending only on $A,B,\alpha$ and $p$ such that for any $t\in [0,T]$, 
$K(t)$ has the interior $\conereg^{\beta/2,C(\omega)}$-property , where $C(\omega)=C_0\|\omega\|_p^{1/2}$. 

In particular, there is a constant 
$$\rho=\frac12 (2C(\omega))^{-2/\beta}=\frac12 (2C_0)^{-2/\beta}\|\omega\|_p^{-1/\beta}$$ 
such that for any $t\in [0,T]$, the set $K(t)$ has the interior 
cone property of parameters $(\rho, 2\rho)$.
\end{Corollary}
%%%%%%%%%%%%

{\bf Proof:} Let us prove the first part of the corollary. Let $K_1$ be such that 
$K_0=K_1+r_0\overline{B}(0,1)$. 
Then $K(t)$ is the reachable set at time $r_0+t$ for the system 
$$
x'(t)=\tilde c(x(t),t)b(t) \qquad |b(t)|\leq 1\;,
$$
starting from $K_1$, where $\tilde c(x,t)=1$ if $t\in [0,r_0]$, and $\tilde c(x,t)=c(x,t-r_0)$
if $t\in (r_0,T+r_0]$ (notice that $\tilde c$ satisfies 
\eqref{BornesVitesse}--\eqref{estiLn}--\eqref{Holder}).
 For this system, Lemma \ref{lem:ccone} shows the result as soon 
as $t\geq C_1C(\omega)^{-2/\beta}.$
Therefore, if we assume that $C_1C(\omega)^{-2/\beta}\leq r_0$, 
which is always possible by increasing $\| \omega \|_p$, 
then the result holds for $K(t)$, for any $t\in [0,T].$

For the second part of the result, let $\theta=2\rho=(2C(\omega))^{-2/\beta}$, $t\in (0,T]$, 
$x\in \partial K(t)$
and $\nu\in \Sc^{N-1}$ be such that $\conereg^{\beta/2,C(\omega)}(x,\nu)\subset K(t)$. 

Since $\theta = (2C(\omega))^{-2/\beta}$, we have $\theta \leq C(\omega)^{-2/\beta}$ and $\rho \leq \theta-C(\omega) \theta^{1+\beta/2}$, so that $\cone_{x,\nu}^{\rho,\theta} \subset \, \conereg^{\beta/2,C(\omega)}.$ This proves that 
the cone $\cone_{x,\nu}^{\rho,\theta}$, with $\theta=2\rho$, is contained in  $K(t)$.
\QED

We now show that the convex set $\conereg^{\delta,C}(\bar x, \nu)$
has a boundary of class
${\mathcal C}^{1,\gamma}$ in a neighborhood of $\bar x$ for some $\gamma>0$. 
Let us fix a frame $\{e_1, \dots, e_N\}$ 
of $\R^N$ such that $\bar x=0$, $\nu = e_N$. We denote by $(x',x_N)$
a generic element of $\R^N$, with $x'\in \R^{N-1}$, $x_N\in \R$. 

%%%%%%%%
\begin{Lemma}\label{lem:Cgraphe} Let $C>0$ and $\delta>0$ be fixed. 
There are constants  $\gamma=\delta/(2+\delta)$, $c=2(2C)^{1/(2+\delta)}$, 
$\tau_0 =(2C)^{-\frac{1}{\delta}}$  and $r_0= (\sqrt3-1)^{\frac{2+\delta}{\delta}}\tau_0$ 
such that the set
$$
\{(x',x_N)\in \R^N\; ;\; |x'|\leq r_0, \; c|x'|^{1+\gamma}\leq x_N\leq \tau_0\}
$$
is contained in $\conereg^{\delta,C}(0, \nu)$.
\end{Lemma}
%%%%%%%

\noindent {\bf Proof:} Note that, by choice of $\tau_0$, 
the map $\tau\to r(\tau)= \tau(1-C\tau^{\delta})$
is nondecreasing on $[0,\tau_0]$. For any $\tau\in (0,\tau_0]$,
the ball $\overline{B}(\tau e_N, r(\tau))$ is contained in  $\conereg^{\delta,C}(0, \nu)$,
which is convex. Let us set $\psi_\tau(x')=\tau-(r^2(\tau)-|x'|^2)^{1/2}$. 
Since the set $\conereg^{\delta,C}(0, \nu)$ is convex, the set
\be\label{InterSet}
\{(x',x_N)\in \R^N\; ; \;|x'|\leq r(\tau), \;  \psi_\tau(x')\leq x_N\leq \tau_0\}
\ee
 is contained in $\conereg^{\delta,C}(0, \nu)$. Indeed, if $|x'|\leq r(\tau ),$
then $(x',\psi_\tau(x'))\in 
\overline{B}(\tau e_N, r(\tau))$ while  $(x',\tau_0)\in \overline{B}(\tau_0 e_N, r(\tau))$. 
 Let $|x'|\leq r_0$ and let us choose $$\tau=(2C)^{-1/(2+\delta)}|x'|^{2/(2+\delta)}.$$ Then
 $\tau\in (0,\tau_0)$ and $|x'|\leq r(\tau)$ (here we use the fact that $|x'|\leq r_0$). 
 Moreover, since $|x'|^2=2C\tau^{2+\delta}$, we get
$$
\begin{array}{rl}
 \psi_\tau(x') \leq & \tau -\left(\tau^2(1-C\tau^{\delta})^2-2C\tau^{2+\delta}\right)^{1/2} \\
\leq & \tau \left[ 1 -\left(1-4C\tau^{\delta}\right)^{1/2} \right]\\
\leq & 2C\tau^{1+\delta}\; = \; (2C)^{1/(2+\delta)}|x'|^{1+\gamma}\;.
\end{array}
$$
Using \eqref{InterSet}, we get that 
any point of the form $(x', x_N)$ with 
$$
|x'|\leq r_0\quad \text{and} \quad c|x'|^{1+\gamma}\leq x_N\leq \tau_0, \quad \text{where} \; 
c=2(2C)^{1/(2+\delta)},
$$ 
belongs to $\conereg^{\delta,C}(0, \nu)$.
\QED

Let us now state a stability property for sets satisfying an interior $\conereg^{\delta,C}$-property: 

%%%%%%%%%%
\begin{Lemma}\label{CvNormGrad} Let $(z_n)$ be a sequence of Lipschitz continuous real-valued maps 
on $\R^N$ which converges uniformly to some $z$. We assume that
$\{z_n\leq 0\}=\{z\leq 0\}$, that
there exist constants $A,B>0$ such that the following inequality holds in the viscosity sense: 
for any $n\in \N$,
$$
\frac{1}{B} \leq |Dz_n(x)|\leq \frac{1}{A}\qquad {\rm in }\; \{0<z_n<T\},
$$
and that there exist $C,\delta>0$ such that for any $x\in \{0<z<T\}$ and any $n$ sufficiently large, 
there is some $\nu\in \Sc^{N-1}$
with $\conereg^{\delta, C}(x,\nu)\subset \{z_n\leq z_n(x)\}$. Then 
$$
\frac{Dz_n(x)}{|Dz_n(x)|}\to \frac{Dz(x)}{|Dz(x)|}\qquad \mbox{\it a.e. in } \{0<z<T\}
$$
and $(|Dz_n| )$ converges to $|Dz|$ in $L^\infty-$weak$-*$ in $\{0<z<T\}$. 
\end{Lemma}
%%%%%%%%%

{\bf Proof: } By standard stability property of viscosity solutions we have that 
$$
\frac{1}{B} \leq |Dz(x)|\leq \frac{1}{A}\qquad {\rm in }\; \{0<z<T\},
$$
in the viscosity and a.e. sense. Note also that, in view of Remark \ref{A consequence}, 
the indicator function of the set $\{0<z_n<T\}$ converges a.e. to the indicator function of $\{0<z<T\}$. 
Let $x$ be such that $z_n$ and $z$ are positive and differentiable at $x$ for any $n$. 
Then $|Dz_n(x)|>0$ for any $n$
and $|Dz(x)|>0$. From the regularity assumption on $z_n$ there exists $\nu_n\in \Sc^{N-1}$
such that $\conereg^{\delta, C}(x,\nu_n)\subset \{z_n\leq z_n(x)\}$. Since $Dz_n(x)$ exists and is nonzero 
and since the set $\conereg^{\delta, C}(x,\nu_n)$ is of class ${\mathcal C}^1$ at $x$ 
(thanks to Lemma \ref{lem:Cgraphe}), one must have
$\nu_n=- Dz_n(x)/|Dz_n(x)|$. Let $\nu$ be the limit of a subsequence of the $(\nu_n)$. Then 
$\conereg^{\delta, C}(x,\nu)\subset \{z\leq z(x)\}$, so that by the same argument as above, 
$\nu=-Dz(x)/|Dz(x)|$. Accordingly any converging subsequence of $Dz_n(x)/|Dz_n(x)|$ 
converges to $Dz(x)/|Dz(x)|$, which shows the 
a.e. convergence of $(Dz_n/|Dz_n|)$ to $Dz/|Dz|$.

Since the $(z_n)$ are uniformly Lipschitz continuous and $(z_n)$ converges uniformly to $z$, $(Dz_n)$
converges  to $Dz$ in $L^\infty-$weak$-*$ in $\{0<z<T\}$. Let $a\in L^1(\R^N, \R^N)$. 
Then  we have on the one hand
$$
\lim_{n\to +\infty} \int_{\{0<z<T\}} \lg a, Dz_n\rg =  \int_{\{0<z<T\}} \lg a, Dz\rg \;.
$$
On the other hand, if we denote by $\xi$ any weak$-*$ limit of a subsequence $(|Dz_{n_k}|)$, we have, from the 
a.e. convergence of $(Dz_n/|Dz_n|)$ to $Dz/|Dz|$, 
$$
\lim_{k\to +\infty} \int_{\{0<z<T\}} \lg a, Dz_{n_k}\rg
= \lim_{k\to +\infty} \int_{\{0<z<T\}} \lg a, \frac{Dz_{n_k}}{|Dz_{n_k}|}\rg |Dz_{n_k}|
= \int_{\{0<z<T\}} \lg a, \frac{Dz}{|Dz|}\rg \xi\;.
$$
This implies that
$$
Dz(x)=\frac{Dz(x)}{|Dz(x)|} \xi(x) \qquad \mbox{\rm a.e. in }\{0<z<T\}\;,
$$
and shows that $\xi=|Dz|$. Hence $(|Dz_n|)$ converges to $|Dz|$ weakly$-*$ in $\{0<z<T\}$. 
\QED

We complete the section by proving that a set with the interior cone property
is the union of a finite number of Lipschitz graphs.

%%%%%%
\begin{Proposition}\label{graphs} Let $(K(t))_{t\in [0,T]}$ be a 
nondecreasing family of compact subsets of $\R^N$, each $K(t)$ having 
the interior cone property
of parameter $(\rho,2\rho)$ for some $\rho>0$. Then for any 
$\bar x\in \R^N$ and any $r\geq \rho$, there is an integer
$C(r,\rho)\leq C(N)r/\rho$ (where $C(N)$ only depends on $N$) 
and, for each $i\in \{1, \dots, C(r,\rho)\}$,
\begin{itemize}
\item a Borel measurable map 
$\Psi_i:B_{N-1}(0,r)\times [0,T]\to \R$, which is 
$\sqrt{15}-$Lipschitz continuous with respect
to the space variable,  
\item and a change of coordinates $O_i:\R^N\to\R^N$ 
(i.e., the composition of a rotation and a translation), with $O_i(0)=\bar x$,
\end{itemize}
such that, for all $t\in [0,T],$
$$
\partial K(t)\cap B(\bar x,r) 
\subset \bigcup_{i=1,\dots,C(r,\rho)} 
\left\{ O_i(x',\Psi_i(x',t))\;,\; x'\in B_{N-1}(0,r)\right\}.
$$
If furthermore the family $(K(t))$ is contained in some ball $\overline{B}(0,M),$ 
then we can take $r=+\infty$ and $C(\rho)\leq C(N)M/\rho$
and we have, for all $t\in [0,T],$
$$
\partial K(t) \subset \bigcup_{i=1,\dots,C(r,\rho)} 
\left\{ O_i(x',\Psi_i(x',t))\;,\; x'\in B_{N-1}(0,M)\right\}.
$$
\end{Proposition}
%%%%%%%

\noindent
An important and straightforward consequence of the fact that
$\partial K(t)$ is piecewise Lipschitz continuous is that
the sets $K(t)$ are of (locally) finite perimeter.\\

{\bf Proof: } We closely follow several arguments of \cite{bclm09}.    
We first observe that if $x \in \partial K$ and 
$\cone_{x,\nu}^{\rho,2\rho} \subset K(t)$, then for all $\nu' \in \Sc^{N-1}$  
verifying $|\nu - \nu'| \leq 1 / 4$, we have 
$\cone_{x,\nu'}^{\rho /  2,2\rho} \subset K(t)$. By compactness of  
$\Sc^{N-1}$, we can cover $\Sc^{N-1}$ with the traces on $\Sc^{N-1}$ of at most $p$  
balls of radius $1 / 4$ centered at $\nu_i$, for some positive constant $p=p(N)$ 
and $1 \leq j \leq p$. Therefore,  
for any $x \in \partial K(t)$, there exists $1\leq j \leq p$ such 
that $\cone_{x,\nu_j}^{\rho / 2,2\rho} \subset K(t)$.
  
\begin{figure}[ht]  
\begin{center}  
\epsfig{file=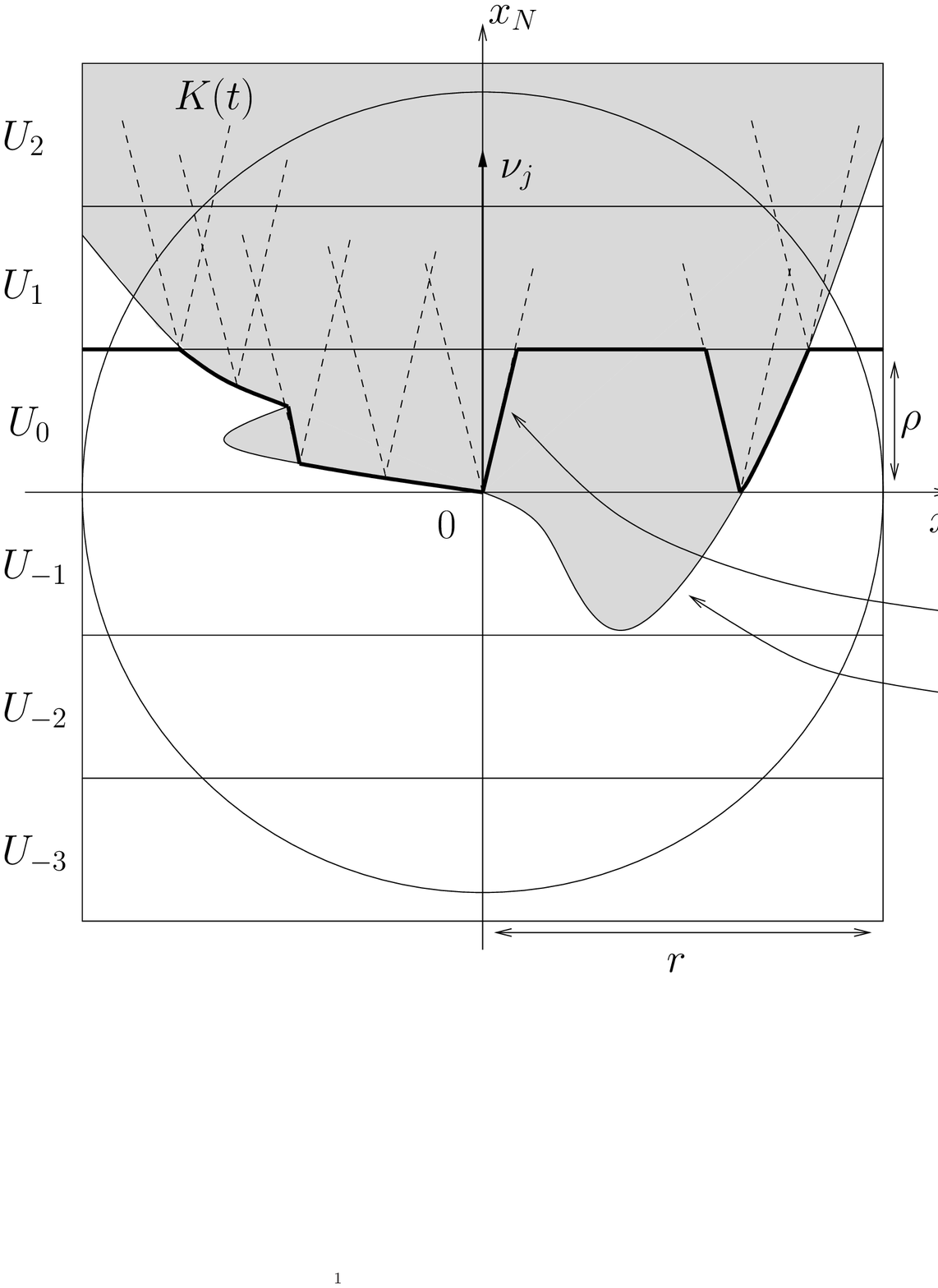, width=8cm}  
\end{center}  
\caption{ \label{fig-decoup}  
{\it
}}  
\end{figure}

Let us now fix $\bar x$ and $1\leq j \leq p$. 
Up to a translation and a rotation of the space, we can assume 
that $\bar x=0$, $\nu_j=(0,\dots,0,1)$.
For any $x\in \R^N$, we write $x=(x',x_N)$ with $x'\in \R^{N-1}$ 
and $x_N\in \R$. For any $t\in [0,T]$ and any 
integer $k$ with $|k|\leq r/\rho+1$, we
set 
$$
U_{k}=B_{N-1}(0,r) \times \left[ k\rho, (k+1)\rho \right]\, ,
$$
$$
A_{j,k}(t)=\left\{ x=(x',x_N)\in \partial K(t)\cap \overline{U_k}
\;; \; \cone_{x,\nu_j}^{\rho / 2,2\rho} \subset K(t)\right\}\, ,
$$
and, for all $y'\in B_{N-1}(0,r),$
$$
\Psi_{j,k}(y',t)= {\rm min}\left\{ 
(k+1)\rho \; , \; 
\mathop{\rm inf}_{x\in A_{j,k}(t)}\psi_x(y')
\right\},
$$
where $\psi_x(y')=\sqrt{15}|y'-x'|+x_N$ is such that
$({\rm graph}\, \psi_x)\cap U_{k}= \cone_{x,\nu_j}^{\rho / 2,2\rho}\cap  U_{k}$
(see Figure \ref{fig-decoup} for an illustration).
%%%%%%%%%%%%%%%%%%%%%%%%%%%%  
% figure: decoupage front en bandes horizontales
%%%%%%%%%%%%%%%%%%%%%%%%%%%%%%%%%%%%%%%%%%%%%%%%%% 
We claim that 
$$
A_{j,k}(t)\cap U_{k}\subset {\rm graph}\, \Psi_{j,k}(\cdot ,t).
$$
Indeed, let $x\in A_{j,k}(t)\cap U_{k}.$ 
If $x\notin {\rm graph}\, \Psi_{j,k}(\cdot ,t),$ then
$\Psi_{j,k}(x' ,t)< \psi_x(x')=x_N.$ Therefore, there exists
$z\in A_{j,k}(t)$ such that $\psi_z(x') <x_N.$
It follows that 
$x\in {\rm int}\, \cone_{z,\nu_j}^{\rho / 2,2\rho}\subset {\rm int}\, K(t)$
and $x$ cannot belong to $\partial K(t),$ which is a contradiction.
This proves the claim. Then we remark that $\Psi_{j,k}(\cdot ,t)$ is a Lipschitz continuous
map with constant $\sqrt{15}$ as the infimum of a family of maps having this property.

This means that $\partial K(t)\cap B(\bar x,r)$ is contained in at most
$p (2r/\rho+2)$ Lipschitz graphs with constant $\sqrt{15}$, which concludes the proof since $r\geq \rho$; indeed this implies that $p (2r/\rho+2) \leq 4p\, r/\rho=:C(r,\rho).$
\QED

%%%%%%%%%%%%%%%%%%%%%%%%%%%%%%%%%%%%%%%%%%%%%%%%%%%%%%%%%%%%%%%%%%%
%%%%%%%%%%%%%%%%%%%%%%%%%
%%%%%%%%%%%%%%%%%%%%%%%%%%%%%%%%%%%%%%%%%%%%%%%%%%%%%%%%%%%%%%%%%%
%%%%%%%%%%%%%%%%%%%%%%%%%%
\section{Representation and a priori estimates for the heat equation}
\label{SecHolderBounds}

The aim of this section is to provide estimates for the following heat equations
\be\label{Heat0}
\left\{\begin{array}{l}
v_t-\Delta v+g(x,t){\mathcal H}^{N-1}\lfloor\Gamma(t)=0\qquad {\rm in }\ \R^N\times (0,T)\ ,\\
v(x,0)=v_0(x)\qquad  {\rm in }\ \R^N\ ,
\end{array}\right.
\ee
and
\be\label{Heat1}
\left\{\begin{array}{l}
v_t-\Delta v+\kappa\bar g(v(x,t)){\mathcal H}^{N-1}\lfloor\Gamma(t)=0\qquad {\rm in }\ \R^N\times (0,T)\ ,\\
v(x,0)=v_0(x)\qquad {\rm in }\ \R^N\ ,
\end{array}\right.
\ee
for a given evolving front $(\Gamma(t))_{t\geq 0}$. 

Throughout the section we work under the following conditions on the data: 
\begin{itemize}

\vspace{.2cm}

\item[{\bf (H1)}] $g:\R^N\times [0,T] \to \R$ is continuous 
and bounded by a constant $M>0$.

\vspace{.2cm}

\item[{\bf (H2)}] $\kappa\in \R$ and $\bar g: \R\to \R$ 
is bounded by $M$ and Lipschitz continuous. 

\vspace{.2cm}

\item[{\bf (H3)}] $v_0$ is Lipschitz continuous and bounded. 

\vspace{.2cm}

\item[{\bf (H4)}] The evolving family $(\Gamma(t))_{t\in [0,T]}$ 
can be represented as
\begin{equation}\label{GammaVSz}
\Gamma(t)=\{x\in \R^N\; ; \; z(x)=t\}\qquad \ {\rm for \ all} \ t\in (0,T)\;.
\end{equation}
where $z:\R^N\to \R$ is Lipschitz continuous and satisfies
\be \label{estiDz}
\frac{1}{B}\leq |Dz(x)|\leq \frac{1}{A} \qquad {\rm in }\; \{0<z<T\}
\ee
in the viscosity sense for some $A,B>0$. 
Furthermore we assume that there is some $\bar \rho>0$ such that the set
$$
K(t)=\{x\in \R^N\; ; \; z(x)\leq t\}
$$
has the interior cone property
of parameter $(\bar \rho,2\bar \rho)$ for all $t\in (0,T)$, and that there exists $M>0$
such that
$$
K(t)\subset \overline{B}(0,M).
$$
\end{itemize}
Let us recall that, thanks to the interior cone condition, $K(t)$ 
is a set of finite perimeter and, moreover, its boundary $\Gamma(t)$ is 
contained in the union of a finite number of Lipschitz graphs (Proposition \ref{graphs}). 

Throughout the section we denote by $C$ a constant which only depends on 
$A,B,N,T,M, \kappa$ and may vary from line to line in the computations. 

%%%%%%%%%%%%%%%%%%%%%%%%%%%%%%%%%%%%%%%%%%
\subsection{Representation and $L^\infty$ bounds for the solution 
of \eqref{Heat0}}

%%%%%%
\begin{Lemma}\label{lem:repsol} There exists a unique solution 
to \eqref{Heat0}. This solution is given,
for all $(x,t)\in \R^N\times [0,T],$ by
\begin{eqnarray*}
v(x,t)=\int_{\R^N}G(x-y,t)v_0(y) \, dy
-\int_0^t\int_{\Gamma(s)}G(x-y,t-s)g(y,s)\, d{\mathcal H}^{N-1}(y)ds,
\end{eqnarray*}
where $G(x,t)=(4\pi t)^{-N/2} e^{-|x|^2/(4t)}$ 
is the kernel of the heat equation, and satisfies the 
uniform bound
\be \label{UniBounds}
|v(x,t)|\leq C(1+|\log(\bar \rho)|)
\qquad \ {for \ all} \ (x,t)\in \R^N\times [0,T],
\ee
where $\bar{\rho}$ is the cone paramater which appears in {\bf (H4)}.
\end{Lemma}
%%%%%%%

{\bf  Proof: } Uniqueness of the solution is clear. 
The term $\int_{\R^N}G(x-y,t)v_0(y) \, dy$ corresponds to the initial datum 
and satisfies the bound
$$
\left| \int_{\R^N}G(x-y,t)v_0(y) \, dy \right| \leq \|v_0\|_\infty\, .
$$
In order to prove the representation formula and the bound for $v$, 
we can therefore assume that $v_0=0.$ Let us set $f_\ep(x,t)={\bf 1}_{K(t)}*G(\cdot,\ep)$ 
(where the convolution is only made with respect to the space variable). 
Then $f_\ep$ is smooth in space and 
strictly converges in the BV sense to ${\bf 1}_{K(t)}$ (see \cite[Def. 3.14]{afp00} 
and \cite[Sect. 5.2]{eg92}). 
In particular, since $\partial K(t)$ is piecewise Lispchitz continuous, the measure
$|Df_\ep(\cdot,t)|dx$ weakly-$*$ converges to ${\mathcal H}^{N-1}\lfloor\Gamma(t)$
(\cite[Prop. 3.62]{afp00}). For all $(x,t)\in \R^N\times [0,T],$ let 
$$
v_\ep(x,t)=-\int_0^t\int_{\R^N} g(y,s)G(x-y,t-s)|Df_\ep(y,s)|dy ds.
$$
Since $|Df_\ep(\cdot,t)|$ is Lipschitz continuous, it is well-known
that $v_\ep$ is a solution of
\be\label{HeatApprox}
(v_\ep)_t-\Delta v_\ep+ g(y,s)|Df_\ep(x,t)|=0\ \  {\rm in}\; \R^N\times(0,T)\;.
\ee
The key step in the proof of \eqref{UniBounds} is the following uniform bound on  $(v_\ep)$:
\be\label{Boundsuep}
|v_\ep(x,t)|\leq C(1+|\log(\bar \rho)|)\ \ {\rm for \ all} \ (x,t)\in \R^N\times [0,T]\;,
\ee
which holds for any $\ep>0$. 
Let us assume for a while that this is true. Then, by the weak-$*$ convergence of $|Df_\ep|dx$ to 
${\mathcal H}^{N-1}\lfloor\Gamma$, $(v_\ep)$ converges pointwise to $v$ in $(\R^N\times(0,T)) \backslash \Gamma$,
hence in $L^1_{loc}(\R^N \times [0,T])$ since it is uniformly bounded in $L^\infty$ thanks to the bound \eqref{Boundsuep}, and $\Gamma$ has zero measure
in $\R^N\times (0,T)$. By \eqref{HeatApprox} $v$ is a solution of \eqref{Heat0}.

\vspace{.2cm}

It remains to prove \eqref{Boundsuep}.  To do this we note that, since $K(t)$ is a set of finite perimeter, we have
$$
|Df_\ep(y,s)|\leq \int_{\Gamma(t)}G(y-x',\ep)d{\mathcal H}^{N-1}(x')\quad
{\rm for \ all} \  (y,s)\in \R^N\times (0,T), \; y\notin \Gamma(s)\;.
$$
Therefore, since $G(x-x',t-s+\ep)=\int_{\R^N}G(x-y,t-s)G(y-x',\ep)dy,$ we get
$$
\begin{aligned}
|v_\ep(x,t)|\; \leq & M \int_0^t\int_{\R^N}\int_{\Gamma(s)} G(x-y,t-s)G(y-x',\ep) \, d{\mathcal H}^{N-1}(x')dy ds\\
\leq &  C\int_0^t\int_{\Gamma(s)} G(x-x',t+\ep-s)d{\mathcal H}^{N-1}(x') ds.
\end{aligned}
$$
Let us split this last integral in two parts, the first one denoted by $I_1$ being the integral between $0$ and $t-\tau$
and the other one, denoted by $I_2$, between $t-\tau$ and $t$ for some $\tau\in (0,t]$. 
Let us first estimate
$$
I_1=C\int_0^{t-\tau} \int_{\Gamma(s)} G(x-y,t+\ep-s)d{\mathcal H}^{N-1}(y)ds\, .
$$
From \eqref{GammaVSz} and Lemma \ref{CoArea} below, we have 
\begin{eqnarray*}
& I_1&=C\int_0^{t-\tau}\int_{\{z=s\}}G(x-y,t+\ep-s)\ d{\mathcal H}^{N-1}(y)ds\\
& &\leq  \frac{C}{A} \left[\int_{\{0<z<t-\tau\}}G(x-y,\ep+\tau)\, dy 
+ \int_0^{t-\tau} \int_{\{0<z<s\}} |G_t(x-y,t+\ep-s)|\ dy ds\right]\;.
\end{eqnarray*}
Note that 
$$
\int_{\{0<z<t-\tau\}}G(x-y,\ep+\tau)\ dy \leq \int_{\R^N}G(x-y,\ep+\tau)\ dy \; =\; 1 \;.
$$
Moreover we have 
\begin{eqnarray*}
&& \int_{\{0<z<s\}}\left|G_t(x-y,t+\ep-s)\right|\, dy \\
&\leq &  
C\int_{\R^N}\left(\frac{1}{(t+\ep-s)^{(N+2)/2}}
+\frac{|y-x|^2}{(t+\ep-s)^{(N+4)/2}}\right) e^{-|y-x|^2/(4(t+\ep-s))} \,dy \\
&\leq & 
C\int_0^{\infty}\left(\frac{r^{N-1}}{(t+\ep-s)^{(N+2)/2}}
+\frac{r^{N+1}}{(t+\ep-s)^{(N+4)/2}}\right) e^{-r^2/(4(t+\ep-s))} \,dr \\
&\leq & 
\frac{C}{t+\ep-s}\int_0^{\infty}(r^{N-1}
+r^{N+1})e^{-r^2} \,dr \; \leq \; \frac{C}{t+\ep-s}\; \leq \; \frac{C}{t-s}.
\end{eqnarray*}
Therefore we get
$$
I_1\; \leq \; C(1+\log(t/\tau))\, .
$$
We now estimate
$$
I_2=C \int_{t-\tau}^t\int_{\Gamma(s)}G(x-y,t+\ep-s)\, d{\mathcal H}^{N-1}(y)ds\, .
$$
From the structure condition on $K(s)$ and Proposition \ref{graphs}, 
there exists an integer $C(\bar \rho)\leq C_1/\bar \rho$ 
(where $C_1$ only depends on $N, M$) and, for each $i\in \{1, \dots, C(\bar \rho)\}$,
\begin{itemize}
\item a Borel measurable map $\Psi_i:B_{N-1}(0,M)\times [0,T]\to \R$, which is $\sqrt{15}-$Lipschitz continuous with respect
to the space variable,  
\item and a change of coordinates $O_i:\R^N\to\R^N$, where $O_i(0)=x,$
\end{itemize}
such that, for all $s\in [0,T],$
$$
\Gamma(s)  \subset \bigcup_{i=1,\dots,C(\bar \rho)} 
\left\{ O_i(y',\Psi_i(y',s)),\; y'\in B_{N-1}(0,M)\right\}.
$$
Therefore, using that 
$$
{\mathcal H}^{N-1}\lfloor \{ (y',\Psi_i(y',s)),y'\in B_{N-1}(0,M)\}
=\sqrt{1+|D\Psi_i(y',s)|^2}\, {\mathcal L}^{N-1}\lfloor B_{N-1}(0,M),
$$ 
we have
$$
I_2 \leq C \sum_{i=1}^{C(\bar \rho)} \int_{t-\tau}^t\int_{B_{N-1}(0,M)}
G((y',\Psi_i(y',s)),t+\ep-s)\sqrt{1+|D\Psi_i(y',s)|^2}\, dy' ds\, .
$$
We deduce that
\begin{eqnarray*}
I_2 & \leq &  
\frac{C}{\bar \rho} \int_{t-\tau}^t\int_{\R^{N-1}}  
\frac{1}{(t+\ep-s)^{N/2}} e^{-|y'|^2/(4(t+\ep-s))}\, dy'ds  \\
&\leq & 
\frac{C}{\bar \rho} \int_{t-\tau}^t\int_0^{+\infty}  
\frac{r^{N-2}}{(t+\ep-s)^{N/2}} e^{-r^2/(4(t+\ep-s))}\, dr ds  \\
&\leq & 
\frac{C}{\bar \rho} \int_{t-\tau}^t\int_0^{+\infty}  
\frac{r^{N-2}}{(t-s)^{1/2}}   e^{-r^2}\, dr ds 
\;\leq  \frac{C\sqrt{\tau}}{\bar\rho}.
\end{eqnarray*}
Putting together the estimates for $I_1$ and $I_2$ gives
$$
|v_\ep(x,t)|\leq C\left(1+\log \left(\frac{t}{\tau}\right)
+\frac{\sqrt{\tau}}{\bar \rho} \right)\, ,
$$
which holds for any $\tau\in (0,t]$. Choosing $\tau=\bar \rho^{2}$ 
if $t\geq  \bar \rho^{2}$ and $\tau=t$ otherwise (in which case 
the decomposition reduces to $I_2$), we finally obtain
\eqref{Boundsuep}.
\QED

The following Lemma, which was used in the proof, is a simple 
consequence of the Coarea formula. 

%%%%%
\begin{Lemma}\label{CoArea} Let $T>0$, $z: \R^N\to \R$ be Lipschitz continuous and such that 
$$
\frac{1}{B}\leq |Dz|\leq  \frac1A  \qquad { a.e. \ in} \ \{0<z<T\}\;.
$$
Let $0\leq s_1<s_2\leq T$ and assume that $\phi: \R^N\times (s_1,s_2) \to \R$ 
is nonnegative and such that $\phi$ and $\phi_t$ are integrable on $\{s_1<z<s_2\}$. 
Then
$$
\begin{aligned}
& \int_{s_1}^{s_2}\int_{\{z=s\}} \phi(x,s)\ d{\mathcal H}^{N-1}(x) ds\\
& \qquad \qquad \leq \frac1A \left[ \int_{\{s_1<z<s_2\}}\phi(x,s_2)\ dx + \int_{s_1}^{s_2} \int_{\{s_1<z<s\}} |\phi_t(x,s)|\ dx ds\right]\;.
\end{aligned}
$$
\end{Lemma}
%%%%%

{\bf Proof : } Let us first assume that $\phi$ is smooth and bounded. 
From the Coarea formula \cite[Sect. 3.4.4]{eg92} we have 
$$
\int_{s_1}^{s_2}\int_{\{z=s\}}\frac{\phi(x,s)}{|Dz(x)|}\,
d{\mathcal H}^{N-1}(x)ds \; =\;  \int_{\{s_1<z<s_2\}}\phi(x,z(x))\, dx
$$
while, by Fubini's Theorem, we get
$$
\begin{aligned}
\int_{s_1}^{s_2} \int_{\{s_1<z<s\}} \phi_t(x,s)\, dxds \; = & \int_{\{s_1<z<s_2\}}\int_{z(x)}^{s_2} \phi_t(x,s)\, dsdx \\
= & \int_{\{s_1<z<s_2\}}\phi(x,s_2)\,dx - \int_{\{s_1<z<s_2\}} \phi(x,z(x))\,dx\; .
\end{aligned}
$$
So 
$$
\begin{aligned}
& \int_{s_1}^{s_2}\int_{\{z=s\}}\frac{\phi(x,s)}{|Dz(x)|}\ d{\mathcal H}^{N-1}(x) ds\\
& \qquad \qquad \leq  \int_{\{s_1<z<s_2\}}\phi(x,s_2)\ dx + \int_{s_1}^{s_2} \int_{\{s_1<z<s\}} |\phi_t(x,s)|\ dx ds\;.
\end{aligned}
$$
Since $|Dz|\leq 1/A$, this gives the result for $\phi$ smooth and bounded. 
The general case follows by regularization.
\QED

We shall need two types of space regularity estimates for the solution $v$ to \eqref{Heat0}. The first one is a continuity estimate with a modulus
$\omega(s)=s(1+|\log(s)|)$: it is required in order to solve unambiguously the eikonal equation with a velocity
$\bar g(v(x,t))$, but is very crude with respect to the $\bar \rho$ dependance; we prove it in Subsection \ref{secmodule}. The second one is merely a H\"{o}lder estimate, but it is much sharper with respect to the $\bar \rho$ dependance. It is the aim of Subsection \ref{secholder}.

%%%%%%%%%%%%%%%%%%%%%%%%%%%%%%%%%%%%%
\subsection{Modulus of continuity in space for the solution of \eqref{Heat0}}\label{secmodule}

%%%%%%%
\begin{Lemma} Let $v$ be the solution of \eqref{Heat0} given by Lemma \ref{lem:repsol}. Then, for any  $x,y\in \R^N$, $t\in [0,T]$, 
\begin{eqnarray}\label{modulelog}
|v(x,t)-v(y,t)|\leq \frac{C}{\bar \rho}|x-y|\left(1+|\log\,|x-y||\right).
\end{eqnarray}
\end{Lemma}
%%%%%%%

{\bf Proof : } We prove the result for $N\geq 3$, the case $N=2$ being similar but simpler. 

The term $x\mapsto \int_{\R^N}G(x-y,t)v_0(y) \, dy$ is Lipschitz continuous with 
constant $\|Dv_0\|_\infty$; we can therefore assume that $v_0=0$ and $t>0$. 

Using again the structure condition on $K(s)$ and Proposition \ref{graphs}, for any $x\in \R^N$,
there is an integer $C(\bar \rho)\leq C_1/\bar \rho$ (where $C_1$ only depends on $N,M$) 
and, for each $i\in \{1, \dots, C(\bar \rho)\}$,
\begin{itemize}
\item a Borel measurable map $\Psi_i:B_{N-1}(0,M)\times [0,T]\to \R$, which is $\sqrt{15}-$Lipschitz continuous with respect
to the space variable,  
\item and a change of coordinates $O_i=R_i\circ\tau_x:\R^N\to\R^N$, where
$\tau_x (z)=z+x,$ $R_i$ is a rotation, such that
$O_i(0)=x$ and
\end{itemize}
$$
\Gamma(s)  \subset \bigcup_{i=1,\dots,C(\bar \rho)} \left\{ O_i(z',\Psi_i(z',s))\;,\; z'\in B_{N-1}(0,M)\right\}\quad {\rm for \ all} \  s\in [0,T]\;.
$$
Setting 
$$
E_i(s)=\left\{ z=(z',\Psi_i(z',s))\;,\; z'\in B_{N-1}(0,M)\right\}
= {\rm graph}(\Psi_i(\cdot,s)_{|B_{N-1}(0,M)}),
$$
for any $h\in\R^N,$ we have
\begin{eqnarray*}
&& |v(x+h,t)-v(x,t)| \\
&\leq&
M \int_{0}^t\int_{\Gamma(s)} 
\left|G(x+h-y,t-s)-G(x-y,t-s)\right|d{\mathcal H}^{N-1}(y)ds\\
&\leq&
C \sum_{i=1}^{C(\bar \rho)}  \int_{0}^t\int_{O_i(E_i(s))}
\left|G(x+h-y,t-s)-G(x-y,t-s)\right|d{\mathcal H}^{N-1}(y)ds\\
&\leq&
C \sum_{i=1}^{C(\bar \rho)}  \int_{0}^t\int_{E_i(s)}
\left|G(x+h-O_i(z),t-s)-G(x-O_i(z),t-s)\right|d{\mathcal H}^{N-1}(z)ds.
\end{eqnarray*}
Let us set $h_i=(h_i',h_{iN}):= R_i^{-1}h$, where $h_i'\in \R^{N-1}$ 
and $h_{iN}\in \R$. 
We note that, for any $z\in \R^N$, $R_i^{-1}(x+h-O_iz)=h_i-z$, so that
$$
G(x+h-O_i(z),t-s)=G(h_i-z,t-s)
$$
because $G(\cdot,t-s)$ has rotational invariance. 
It follows that
\begin{eqnarray*}
&& |v(x+h,t)-v(x,t)| \\
&\leq&
C \sum_{i=1}^{C(\bar \rho)}  \int_{0}^t\int_{E_i(s)}
\left|G(h_i-z,t-s)-G(-z,t-s)\right|\, d{\mathcal H}^{N-1}(z)ds\\
&\leq&
C \sum_{i=1}^{C(\bar \rho)}  \int_{0}^t\int_{B_{N-1}(0,M)}
\left|G((h_i'-z',h_{iN}-\Psi_i(z',s)),t-s)\right.\\
&& \hspace*{3cm} \left. -G((-z',-\Psi_i(z',s)),t-s)\right|
\sqrt{1+|D\Psi_i(z',s)|^2}\,dz' ds
\end{eqnarray*}
since ${\mathcal H}^{N-1}\lfloor E_i(s) 
= \sqrt{1+|D\Psi_i(y',s)|^2}{\mathcal L}^{N-1}\lfloor B_{N-1}(0,M).$

\vspace{.2cm}

We recall that $|D\Psi_i(z',s)|\leq \sqrt{15}$ and introduce
$$
D_i(s)=\bigcup_{\sigma\in [0,1]} B_{N-1}(\sigma h_i',|h|(t-s)^{1/4})
$$ 
in order to split the latter integral into two parts. We get
\begin{eqnarray*}
&&|v(x+h,t)-v(x,t)| \\
& \leq &  C\sum_{i=1}^{C(\bar \rho)}
\int_{0}^t\int_{D_i(s)}\left|G((h_i'-z',h_{iN}-\Psi_i(z',s)),t-s)\right.\\
&& \hspace{5cm} \left. -G((-z',-\Psi_i(z',s)),t-s)\right| \,dz' ds\\
&& +  |h|\int_0^1\int_{0}^t\int_{\R^{N-1}\backslash D_i(s)}
\left|DG((\sigma h_i'-z',\sigma h_{iN}-\Psi_i(z',s)),t-s)|\, dz'dsd\sigma \right. \\
&=&  C\sum_{i=1}^{C(\bar \rho)} \left( I_i+|h| \, J_i\right)\, .
\end{eqnarray*}

Let us fix $i\in \{1, \dots, C(\bar \rho)\}$ and estimate $I_i$. 
Without loss of generality we can assume that $h_i$ 
belongs to the plane spanned by $e_1$ and $e_N$. 
Then,
$$
D_i(s) \subset \R\times B_{N-2}(0,|h|(t-s)^{1/4})\;,
$$
and setting $z'=(z_1,z'')$ with $z_1 \in\R,$ $z''\in \R^{N-2}$, we have
\begin{eqnarray*}
I_i&  \leq &
C\int_0^t\int_{\R}\int_{ B_{N-2}(0,|h|(t-s)^{1/4}) }
\frac{1}{(t-s)^{N/2}} 
e^{-\frac{|h_1-z_1|^2+|h''-z''|^2+|h_{iN}-\Psi_i(z',s)|^2}{4(t-s)}} \, dz'' dz_1 ds\\
&& +
C\int_0^t\int_{\R}\int_{ B_{N-2}(0,|h|(t-s)^{1/4}) }
\frac{1}{(t-s)^{N/2}} 
e^{-\frac{|z_1|^2+|z''|^2+|\Psi_i(z',s)|^2}{4(t-s)}}\, dz'' dz_1  ds\\
&\leq&
 C\int_0^t\int_{0}^{|h|(t-s)^{1/4}} 
\left( \int_\R \left(e^{-\frac{|h_1-z_1|^2}{4(t-s)}}+ e^{-\frac{|z_1|^2}{4(t-s)}}\right) \, dz_1\right)
\frac{r^{N-3}}{(t-s)^{N/2}}
e^{-\frac{r^2}{4(t-s)}}\, dr ds\\
&\leq &
C\int_0^t\int_{0}^{|h|(t-s)^{-1/4}} \frac{r^{N-3}}{(t-s)^{1/2}}
e^{-r^2/4} \, dr ds\\
&\leq &
 C\int_{0}^{+\infty} \int_{0\vee (t-(|h|/r)^4)}^t  \frac{r^{N-3}}{(t-s)^{1/2}}
e^{-r^2/4} \, dsdr\\
&\leq &
 C\left( \int_{0}^{|h|/t^{1/4}} r^{N-3}e^{-r^2/4}\, t^{1/2}\, dr
+ \int_{|h|/t^{1/4}}^{+\infty}r^{N-5}e^{-r^2/4} |h|^2 \, dr\right)\; .
\end{eqnarray*}
Let $M_N=\mathop{\rm sup}_{[0,+\infty)}r^{N-3}e^{-r^2/4}$ (recall that
$N\geq 3$ by assumption). Then
\begin{eqnarray}
I_i&  \leq &
CM_N  \left(|h|\, t^{1/4}
+ \int_{|h|/t^{1/4}}^{+\infty} \frac{ |h|^2}{r^2}\, dr
\right)\nonumber\\
&\leq &
CM_N T^{1/4}  |h|\nonumber\\
&=&
C |h|.\label{estiI}
 \end{eqnarray}

We now estimate $J_i$. We have
$$
\begin{aligned}
& \left|DG((\sigma h_i'-z',\sigma h_{iN}-\Psi_i(z',s)),t-s)\right|\\
\leq \; &
C \frac{|\sigma h_i'-z'|+|\sigma h_{iN}-\Psi_i(z',s)|}{(t-s)^{(N+2)/2}}
e^{-|\sigma h_i'-z'|^2/(4(t-s))}
e^{-|\sigma h_{iN}-\Psi_i(z',s)|^2/(4(t-s))},
\end{aligned}
$$
with
$$
|\sigma h_{iN}-\Psi_i(z',s)|e^{- |\sigma h_{iN}-\Psi_i(z',s)|^2/(4(t-s))}
\leq C(t-s)^{1/2}\;.
$$
Since $\R^N\setminus D_i(s)\subset \R^N\setminus B_{N-1}(0,|h|(t-s)^{1/4}),$
we get
\begin{equation}
\begin{aligned}
J_i &\leq C\int_0^1\int_{0}^t\int_{|h|(t-s)^{1/4}}^{+\infty}
\left(\frac{r^{N-1}}{(t-s)^{(N+2)/2}}
+\frac{r^{N-2}}{(t-s)^{(N+1)/2}}\right)e^{-r^2/(4(t-s))}\,dr  dsd\sigma \\
& \leq C\int_{|h|t^{-1/4}}^{+\infty}\int_{0}^{t-(|h|/r)^4}
\frac{r^{N-1}+r^{N-2}}{t-s} e^{-r^2/4} \,ds dr \\
& \leq C \int_{|h|t^{-1/4}}^{+\infty}(r^{N-1}+r^{N-2})
\log \left( \frac{tr^4}{|h|^4} \right) e^{-r^2/4}\,dr\\
& \leq C \int_{0}^{+\infty}(r^{N-1}+r^{N-2})
\left( |\log(T)|+|\log(r)|+|\log(|h|) \right) e^{-r^2/4}\,dr\\
& \leq C(1+ |\log |h||).\label{estiJcastwo}
\end{aligned}
\ee
Finally, combining \eqref{estiI}, \eqref{estiJcastwo} and the bound $C(\bar \rho)\leq C_1/\bar\rho$, we obtain
\eqref{modulelog}.
\QED

%%%%%%%%%%%%%%%%%%%%%%%%%%%%%%%%%%%%%
%%%%%%%%%%%%%%%%%%%%%%%%%%%%%%%%%%%%%
\subsection{H\"{o}lder estimate for the solution of \eqref{Heat0}}\label{secholder}

\begin{Lemma}[H\"{o}lder bounds] \label{lem:holder} Let $v$ be the solution of \eqref{Heat0} given by Lemma \ref{lem:repsol}. Then, for any $t\in [0,T]$, $x,y\in \R^N$, 
\begin{equation}\label{eq:estiH}
|v(x,t)-v(y,t)|\leq C(1+|\log(\bar \rho)|)(\bar \rho)^{-\frac14}\, |x-y|^{\frac12}\,.
\end{equation}
\end{Lemma}

{\bf Proof: } The main part of the proof consists in showing the following local H\"{o}lder inequality: 
for any $t\in [0,T]$, $x,h\in \R^N$ with $|h|\leq \sqrt{\bar \rho}/4$, we have
$$
|v(x+h,t)-v(x,t)|\leq C(\bar \rho)^{-\frac14}\, |h|^{\frac12}\,.
$$
We will complete the proof of \eqref{eq:estiH} by using Lemma \ref{lem:repsol}.
 
The term $x\mapsto \int_{\R^N}G(x-y,t)v_0(y) \, dy$ is Lipschitz continuous 
with constant $\|Dv_0\|_\infty$, and therefore locally $1/2$-H\"{o}lder 
continuous; we can assume that $v_0=0$ and $t>0$. Then
\begin{eqnarray*}
& & |v(x+h,t)-v(x,t)| \\
&\leq & M \left[\; |h|\int_0^1\int_0^{t-\tau}\int_{\Gamma(s)}
\left|DG(x+\sigma h-y,t-s)\right|d{\mathcal H}^{N-1}(y)\, dsd\sigma \right. \\
&& +\int_{t-\tau}^t\int_{\Gamma(s)\backslash B(x,\bar r)}
\left(G(x-y,t-s)+G(x+h-y,t-s)\right)d{\mathcal H}^{N-1}(y)\, dsd\sigma \\
&& \left. + \int_{t-\tau}^t\int_{\Gamma(s)\cap B(x,\bar r)}
\left(G(x-y,t-s)+G(x+h-y,t-s)\right)d{\mathcal H}^{N-1}(y)\, 
dsd\sigma \;\right]\\
& = & \|g\|_\infty \left[ \; |h|J_1+J_2+J_3\; \right]
\end{eqnarray*}
where $\bar r,\tau>0$ are chosen such that
$$
\bar r= \sqrt{\bar \rho}\qquad {\rm and }\qquad \tau = |h| \sqrt{\bar \rho} \;.
$$
Since $|h|\leq \sqrt{\bar \rho}/4$, we have $\tau\leq \bar \rho/4$ 
and $\bar r/\sqrt{\tau}\geq 2$. If $\tau >t$, the decomposition 
reduces to $J_2+J_3$ with $\tau=t.$

\vspace{.2cm}

In order to estimate $J_1$, we argue as for $I_1$ in the proof of the estimate \eqref{Boundsuep}: we have
$$
\begin{aligned}
\int_0^1\int_0^{t-\tau}&\int_{\Gamma(s)}\left|DG(x+\sigma h-y,t-s)\right|
d{\mathcal H}^{N-1}(y)\, dsd\sigma  \\
&  \leq\; C\int_0^1\int_0^{t-\tau}\int_{\Gamma(s)}
\frac{|y-x-\sigma h|}{(t-s)^{(N+2)/2}}
e^{-|y-x-\sigma h|^2/(4(t-s))} d{\mathcal H}^{N-1}(y)\, dsd\sigma,
\end{aligned}
$$
where, using Lemma \ref{CoArea}, we have for any $\sigma\in (0,1)$:
\begin{eqnarray*}
&& \int_0^{t-\tau}\int_{\Gamma(s)}\frac{|y-x-\sigma h|}{(t-s)^{(N+2)/2}}
e^{-|y-x-\sigma h|^2/(4(t-s))} \, d{\mathcal H}^{N-1}(y)ds \\
& \leq & \frac1A \left[ 
\int_{K(t-\tau )}\frac{|y-x-\sigma h|}{\tau^{(N+2)/2}}
e^{-|y-x-\sigma h|^2/(4\tau)} \,dy \right.\\
& & \left.+  \; C  \int_0^{t-\tau}\int_{K(s)}
\left(\frac{|y-x-\sigma h|}{(t-s)^{(N+4)/2}}
+\frac{|y-x-\sigma h|^3}{(t-s)^{(N+6)/2}}\right)
e^{-|y-x-\sigma h|^2/(4(t-s))}\, dyds\right].
\end{eqnarray*}
Since, for any $\sigma\in (0,1)$, we have
$$
\int_{K(t-\tau)}\frac{|y-x-\sigma h|}{\tau^{(N+2)/2}}
e^{-|y-x-\sigma h|^2/(4\tau)} \,dy 
\leq   
\int_0^{+\infty}\frac{r^N}{\tau^{(N+2)/2}}e^{-r^2/(4\tau)}dr
 \leq C\tau^{-\frac12}
$$
and
\begin{eqnarray*}
&& \int_0^{t-\tau}\int_{K(s)}\left(\frac{|y-x-\sigma h|}{(t-s)^{(N+4)/2}}
+\frac{|y-x-\sigma h|^3}{(t-s)^{(N+6)/2}}\right)
e^{-|y-x-\sigma h|^2/(4(t-s))} dyds \\
& \leq & C \int_0^{t-\tau}\int_0^{+\infty}\frac{r^N+r^{N+2}}{(t-s)^{3/2}}
e^{-r^2/4}\, dr ds\\
& \leq & C\ \tau^{-1/2},
\end{eqnarray*}
we get
$$
J_1\leq C\, \tau^{-1/2}\,.
$$

For $J_2$ we use the same strategy of proof: from Lemma \ref{CoArea} we have, for any $\epsilon\in (0,\tau)$, 
\begin{eqnarray*}
&& \int_{t-\tau}^{t-\epsilon}\int_{\Gamma(s)\backslash B(x,\bar r)}
\left(G(x-y,t-s)+G(x+h-y,t-s)\right)d{\mathcal H}^{N-1}(y)\, dsd\sigma \\
& \leq & \frac1A \left[ \int_{\{t-\tau<z<t-\epsilon\}} 
{\bf 1}_{\R^N\backslash B(x,\bar r)}(y)
(G(x-y,\epsilon)+G(x+h-y,\epsilon))\, dy\right. \\
& & \left.  
+ \int_{t-\tau}^{t-\epsilon}\int_{\{t-\tau<z<s\}}
{\bf 1}_{\R^N\backslash B(x,\bar r)}(y) 
\left| G_t(x-y,t-s)+G_t(x+h-y,t-s)\right|\, dy ds\right].
\end{eqnarray*}
It is easily seen that 
$$
\lim_{\epsilon\to 0} \int_{\{t-\tau<z<t-\epsilon\}} 
{\bf 1}_{\R^N\backslash B(x,\bar r)}(y)
(G(x-y,\epsilon)+G(x+h-y,\epsilon))\, dy =0,
$$
because $\bar r$ is larger than $4|h|$. On the other hand
$$
\begin{aligned}
&\int_{t-\tau}^{t-\epsilon}\int_{\{t-\tau<z<s\}}
{\bf 1}_{\R^N\backslash B(x,\bar r)}(y) 
\left|G_t(x-y,t-s)+G_t(x+h-y,t-s)\right|\, dy  ds\\
& \leq \; C \int_{t-\tau}^{t} \int_{\bar r/(2(t-s)^{1/2})}^{+\infty} 
\frac{r^{N-1}+r^{N+1}}{t-s}e^{-r^2/4} \,dr ds \\
& \leq \; C \int_{\bar r/(2\sqrt{\tau})}^{+\infty}
\int_{t-\tau}^{t-{\bar r}^2/(4r^2)} 
\frac{r^{N-1}+r^{N+1}}{t-s}e^{-r^2/4}\, ds dr \\
& \leq \; C \int_{\bar r/(2\sqrt{\tau})}^{+\infty}(r^{N-1}+r^{N+1})
\log\left(\frac{4\tau r^2}{{\bar r}^2}\right)
e^{-r^2/4} \, dr\\
& \leq \; C \frac{\sqrt{\tau}}{\bar r}\int_{1}^{+\infty}(r^{N}+r^{N+2})
\log ( r^2)e^{-r^2/4} \, dr
\end{aligned}
$$
because $\bar r/\sqrt{\tau}$ is larger than 2. 
So $J_2\leq C\sqrt{\tau}/\bar r$. 
\vspace{.2cm}

In order to estimate $J_3$ we use the structure of $K(s)$: 
from Proposition \ref{graphs}, 
there exists an integer $C(\bar r,\bar \rho)\leq C_1\bar r/\bar \rho$ 
(where $C_1$ only depends on $N$) and, for each 
$i\in \{1, \dots, C(\bar r,\bar \rho)\}$,
\begin{itemize}
\item a Borel measurable map $\Psi_i:B_{N-1}(0,\bar r)\times [0,T]\to \R$, 
which is Lipschitz continuous with constant $\sqrt{15}$ with respect
to the space variable,  
\item and a change of coordinates $O_i=R_i\circ \tau_i:\R^N\to\R^N$ 
(where $R_i$ is a rotation and $\tau_x$ is a translation),
with $O_i(0)=x,$
\end{itemize}
 such that, for all $s\in [0,T],$
$$
\Gamma(s)\cap B(x,\bar r)  \subset \bigcup_{i=1,\dots,C(\bar r,\bar \rho)} 
\left\{ O_i(z',\Psi_i(z',t))\;,\; z'\in B_{N-1}(0,\bar r)\right\}.
$$ 
Let us set, for any $i\in \{1, \dots, C(\bar r, \bar \rho)\}$, 
$h_i=(h_i',h_{iN}):=R_i^{-1}h$ where $h_i'\in \R^{N-1}$ and $h_{iN}\in \R$. 
Then 
$$
\begin{aligned}
J_{3}&\leq\; \sum_{i=1}^{C(\bar r,\bar \rho)} \int_{t-\tau}^t
\int_{B_{N-1}(0,\bar r)} [G((-z',\Psi_i(z',s)),t-s)\\
& \hspace{1.5cm} +G((h_i'-z',h_{iN}-\Psi_i(z',s)),t-s)]
\sqrt{1+|D\Psi_i(z',s)|^2} \, dz'ds \\
&= \; \sum_{i=1}^{C(\bar r,\bar \rho)} J_{3,i}\;.
\end{aligned}
$$
Let us fix $i\in \{1,\dots, C(\bar r,\bar \rho)\}$. 
Since $|h|\leq \sqrt{\bar \rho}/4= \bar r/4$,  we have
$$
\begin{aligned}
J_{3,i} & \leq \; C \int_{t-\tau}^t \int_{B_{N-1}(0,\bar r)} 
\frac{e^{-|z'|^2/(4(t-s))}+ 
e^{-|h_i'-z'|^2/(4(t-s))}}{(t-s)^{N/2}}\, dz'ds \\
& \leq \; C \int_{t-\tau}^t \int_{B_{N-1}(0,2\bar r)} 
\frac{e^{-|z'|^2/(4(t-s))}}{(t-s)^{N/2}}\, dz'ds\; . \\
\end{aligned}
$$
It follows that
$$
\begin{aligned}
J_{3,i}&  \leq \;  C\int_{t-\tau}^t\int_{0}^{2{\bar r}/(t-s)^{1/2}}  
\frac{r^{N-2}}{(t-s)^{1/2}} e^{-r^2/4} \,dr ds \\
& \leq \; C\int_{0}^{+\infty}\int_{(t-\tau)\vee(t-(2\bar r)^2/r^2)}^{t}  
\frac{r^{N-2}}{(t-s)^{1/2}} e^{-r^2/4} \,ds dr  \\
& \leq \; C\left(\sqrt{\tau} \int_{0}^{2\bar r/\sqrt{\tau}} 
r^{N-2}e^{-r^2/4}\, ds dr
+2\bar r\int_{2\bar r/\sqrt{\tau}}^{+\infty} 
r^{N-3}e^{-r^2/4} \,ds dr\right)  \\
& \leq \; C\left(\sqrt{\tau}\int_{0}^{+\infty} 
r^{N-2}e^{-r^2/4}\,ds dr+\sqrt{\tau}\int_{4}^{+\infty} 
r^{N-2}e^{-r^2/4} \,ds dr\right)\, 
\end{aligned}
$$
since $\bar r/\sqrt{\tau}\geq 2$. 
Accordingly
$$
J_3\leq C \frac{\bar r}{\bar \rho} \, \sqrt{\tau}.
$$
Therefore
$$
|h|J_1+J_2+J_3\leq C\left( \frac{|h|}{\sqrt{\tau}} 
+ \frac{\sqrt{\tau}}{\bar r}+ \frac{\bar r \sqrt{\tau}}{\bar \rho}\right).
$$
With the choice of $\bar r= \sqrt{\bar \rho}$ and $\tau = |h| \sqrt{\bar \rho}$ we get
\begin{equation}\label{titi1}
|v(x+h,t)-v(x,t)|\leq C(\bar \rho)^{-\frac14}\, |h|^{\frac12}
\ \ {\rm for \ all} \
(h,t)\in \R^N\times [0,T] \ {\rm with} \ |h|\leq \sqrt{\bar \rho}/4.
\end{equation}
Now recall that, according to Lemma \ref{lem:repsol}, we have
\begin{equation}\label{titi2}
|v(x,t)|\leq C(1+|\log(\bar \rho)|)
\ \ {\rm for \ all} \ (x,t)\in \R^N\times (0,T).
\end{equation}
Combining \eqref{titi1} and \eqref{titi2} then implies \eqref{eq:estiH}. 
\QED

%%%%%%%%%%%%%%%%%%%%%%%%%%%%%%%%%%%%%
\subsection{Existence, bounds and H\"{o}lder estimate for the solution of \eqref{Heat1}}

%%%%%%%%
\begin{Lemma}\label{ExStabHeat1} Equation \eqref{Heat1} has a 
unique solution $v:\R^N\times [0,T]\to\R$, 
given by
$$
v(x,t)=\int_{\R^N}G(x-y,t)v_0(y) dy
-\kappa\int_0^t\int_{\Gamma(s)}G(x-y,t-s)\bar g(v(y,s)) 
d{\mathcal H}^{N-1}(y)ds.
$$
For all $x,y\in \R^N,$ $t,s\in [0,T],$ $v$ satisfies the following
estimates.

\vspace{.2cm}

$(i)$ Uniform $L^\infty$ bound:
\begin{eqnarray}
 |v(x,t)|\leq C(1+|\log(\bar \rho)|),
\label{BoundsHeat1}
\end{eqnarray}

\vspace{.2cm}

$(ii)$ Space modulus of continuity:
\begin{eqnarray}
|v(x,t)-v(y,t)|\leq \frac{C}{\bar \rho} \, |x-y|(1+|\log |x-y||),
\label{slogsHeat1}
\end{eqnarray}

\vspace{.2cm}

$(iii)$ Space-time H\"older continuity:
\begin{eqnarray}
&& |v(x,t)-v(y,t)|\leq C(1+|\log(\bar \rho)|)\, 
(\bar \rho)^{-1/4}\, |x-y|^{1/2}, 
\label{HolderHeat1}\\
&& |v(x,t)-v(x,s)|\leq \frac{C}{\bar \rho}\, (1+|\log |h| |)\, |t-s|^{1/2}.
\label{THolderHeat1}
\end{eqnarray}
\end{Lemma}
%%%%%%%%

{\bf Proof: } The existence, uniqueness, representation and space 
estimates for the solution of \eqref{Heat1} follow from Banach fixed 
point theorem and Lemmata \ref{lem:repsol}--\ref{lem:holder}. 

Let us now check the time estimate; we fix $0 \leq s\leq t \leq T$ and set $h=t-s$. We note that, from the uniqueness of the solution we have,
for any $x\in \R^N$, 
$$
\begin{aligned}
v(x,t+h) = \int_{\R^N} &G(x-y,h)v(y,t)dy\\ &-\kappa\int_0^h\int_{\Gamma(t+s)}G(x-y,h-s)\, \bar g(v(y, t+s))d{\mathcal H}^{N-1}(y)ds\;.
\end{aligned}
$$
Since $v$ satisfies \eqref{slogsHeat1}, we get from standard estimates 
on the heat flow  that
$$
\left|\int_{\R^N} G(x-y,h)v(y,t)dy-v(x,t)\right|
\leq \frac{C}{\bar \rho} (1+|\log |h| |)h^{\frac12}.
$$
From the structure condition on $K(s)$ and Proposition \ref{graphs}
(see the computations in the proof of Lemma \ref{lem:repsol} for details),
we have
\begin{eqnarray*}
&&\left|\int_0^h\int_{\Gamma(t+s)}G(x-y,h-s)\, \bar g(v(y, t+s)) 
\, d{\mathcal H}^{N-1}(y)ds\right|\\
& \leq & \frac{C}{\bar \rho} \int_0^h\int_{\R^{N-1}}
\frac{1}{(h-s)^{N/2}}e^{-|y'-x'|^2/(4(h-s))} \,dy'ds\\
%& \leq &  \frac{C}{\bar \rho} \int_0^h\int_{0}^{+\infty}
%\frac{\rho^{N-2}}{(h-s)^{N/2}}e^{-\rho^2/(2(h-s))}\,d\rho ds\\
& \leq &  \frac{C}{\bar \rho} \int_0^h\int_{0}^{+\infty}
\frac{r^{N-2}}{(h-s)^{1/2}}e^{-r^2/4} \,dr ds\\
& \leq & \frac{C\sqrt{h}}{\bar \rho}.
\end{eqnarray*}
Putting together the above estimates gives \eqref{THolderHeat1}. 
\QED

%%%%%%%%%%%%%%%%%%%%%%%%%%%%%%%%%%%%%%%%%%%%%%%%%%%%%%%%%%%%%%%%%%%%%%%%%%%%%%%%%%%%%%%%%%%%%%%%%%%%
%%%%%%%%%%%%%%%%%%%%%%%%%%%%%%%%%%%%%%%%%%%%%%%%%%%%%%%%%%%%%%%%%%%%%%%%%%%%%%%%%%%%%%%%%%%%%%%%%%%%
\section{Stability and existence of solutions for the system \eqref{intro:FPP}}

We start with an {\it a priori} stability property for the solution and then prove our main result.

%%%%%%%%%%%%%%%%%%%%%%%%%%%%%%%%%%%%%%%%%%%%%%%%%%
\subsection{A stability property}

We first investigate the convergence of the solution of
$$
\left\{\begin{array}{ll}
(u_n)_t=c_n(x,t)|Du_n| & {\rm in }\; \R^N\times (0,T)\\
(v_n)_t-\Delta v_n+\kappa\bar g(v_n){\mathcal H}^{N-1}\lfloor\{u_n(\cdot,t)=0\}=0 & {\rm in }\; \R^N\times (0,T)\\
v_n(x,0)=v_0(x), \; u_n(x,0)=u_0(x) & {\rm in }\; \R^N
\end{array}\right.
$$
to the solution of
$$
\left\{\begin{array}{ll}
u_t=c(x,t)|Du| & {\rm in }\; \R^N\times (0,T)\\
v_t-\Delta v+\kappa\bar g(v){\mathcal H}^{N-1}\lfloor\{u(\cdot,t)=0\}=0 & {\rm in }\; \R^N\times (0,T)\\
v(x,0)=v_0(x), \; u(x,0)=u_0(x) & {\rm in }\; \R^N
\end{array}\right.
$$
as $(c_n)$ converges to $c$. 

%%%%%
\begin{Lemma}\label{continuite} 
Let us assume that
\begin{itemize}
\item For any $n\in \N$, the velocity $c_n:\R^N\to [0,T]$ satisfies \eqref{BornesVitesse}--\eqref{estiLn}--\eqref{Holder} with fixed $\alpha>1/p$ and 
modulus $\omega$.

\item The sequence $(c_n)$ converges a.e. to some $c:\R^N\times [0,T] \to\R$. 
\end{itemize}

Then $(v_n)$ converges locally uniformly to $v$ in $\R^N\times[0,T]$.
\end{Lemma}
%%%%%

{\bf Proof: } Without loss of generality we can assume that $v_0=0$.
Let us set as usual
$$
K_n(t)=\{u_n(\cdot,t)\geq 0\}, \
\Gamma_n(t)=\{u_n(\cdot,t)= 0\}, \ 
z_n(x)=\inf\{t\geq 0\; ; \; x\in K_n(t)\},
$$
and 
$$
K(t)=\{u(\cdot,t)\geq 0\},\ 
\Gamma(t)=\{u(\cdot,t)= 0\},\ 
z(x)=\inf\{t\geq 0 \; ; \; x\in K(t)\}.
$$
From Proposition \ref{existence} we know that $(u_n)$ converges locally 
uniformly to $u$. 

We claim that this implies that  
$(z_n)$ converges uniformly to $z$ in $\{0<z<t\}$.
Indeed, $u_n(x,z_n(x))=0$ for all $n$ and, passing to the
limit, we get $u(x, {\rm lim\,inf}\,z_n(x))=0.$
Thus ${\rm lim\,inf}\,z_n(x)\geq z(x).$
Now, let $x\in \{0<z<t\}.$ From Proposition~\ref{zLipschitz},
for every $\epsilon,$ there exists $x_\epsilon$ such that 
$|x-x_\epsilon|<\epsilon$ and $u(x_\epsilon,z(x))>0.$
For $n$ sufficiently large, we also have
$u_n(x_\epsilon,z(x))>0$ and therefore $z_n(x_\epsilon)<z(x).$
It follows that $ {\rm lim\,sup}\,z_n(x_\epsilon) \leq z(x).$
Applying again Proposition~\ref{zLipschitz},
we get $-|x-x_\epsilon|/A+  {\rm lim\,sup}\,z_n(x)\leq z(x).$
We conclude by sending $\epsilon$ to 0.

Corollary \ref{cor:ccone} 
states that there is some $\bar \rho>0$ such that each $K_n(t)$ has 
the interior cone property of parameter $(\bar \rho,2\bar \rho)$
and that, for any $x\in \partial K_n(t)$, there is a vector $\nu\in \R^N$ such that $|\nu|=1$ and the set 
$\conereg^{\beta/2,C}(x, \nu)$ is contained in $K_n(t)$, where 
$C=C_0\|\omega\|_p^{1/2}$ and $\beta=\alpha-1/p$. 
Then Lemma \ref{CvNormGrad} implies that $|Dz_n|$ weakly-$*$ 
converges to $|Dz|$ in $\{ 0<z<T \}$. 

By the representation formula for the solution of \eqref{Heat1} 
(Lemma \ref{ExStabHeat1})
and Lemma \ref{extremality} $(2)$,
$$
v_n(x,t) \; = \; -\kappa\int_0^t\int_{\{z_n=s\}}G(x-y,t-s)\bar g(v_n(y,s))d{\mathcal H}^{N-1}(y)ds\, .
$$
From the estimates of Lemma \ref{ExStabHeat1} we know that the $v_n$ 
are uniformly bounded and uniformly H\"{o}lder continuous.
So, up to some subsequence, we can assume that $(v_n)$ uniformly 
converges to some $\bar v$. Our aim is to show that $\bar v=v$. 

Fix $x\in \R^N$ and let $\theta\in (0,t)$ be small. 
Then, following for instance the estimates 
obtained for the proof of \eqref{THolderHeat1}, one easily checks that
$$
\begin{aligned}
&\left|v_n(x,t) +\kappa \int_{0}^{t-\theta}\int_{\{z_n=s\}}G(x-y,t-s)
\bar  g(v_n(y,s))\, d{\mathcal H}^{N-1}(y)ds\right| \\
& \quad \leq \; |\kappa|\|\bar g\|_\infty\int_{t-\theta}^{t}
\int_{\{z_n=s\}}G(x-y,t-s)\, d{\mathcal H}^{N-1}(y)ds\\ 
&\quad \leq \; C(\bar \rho) \, \theta^{1/2}\; .
\end{aligned}
$$
By the Coarea formula, we have
$$
\begin{aligned}
\int_{0}^{t-\theta}&\int_{\{z_n=s\}}G(x-y,t-s)\, \bar g(v_n(y,s)) 
d{\mathcal H}^{N-1}(y)\, ds\\
& = \int_{\{0<z_n<t-\theta\}}G(x-y,t-z_n(y))\, 
\bar g(v_n(y,z_n(y)))|Dz_n(y)| dy.
\end{aligned}
$$
In this expression, 
$$
G(x-\cdot,t-z_n(\cdot)) \, \bar g(v_n(\cdot,z_n(\cdot))) 
\underset{n\to +\infty}{\longrightarrow} G(x-\cdot,t-z(\cdot)) \, 
\bar g( \bar v(\cdot,z(\cdot)))
$$ 
uniformly in $\{0<z<t-\theta\}$
while $(|Dz_n|)$ converges weakly-$*$ to $|Dz|$. Moreover, 
by Remark \ref{A consequence}, the front $\Gamma(s)$ has zero 
measure for any $s$. Therefore, the indicator function of $\{0<z_n<t-\theta\}$
converges to the indicator function of $\{0<z<t-\theta\}$ almost everywhere. 
It follows that
$$
\begin{aligned}
&\underset{n\to +\infty}{\lim} \int_{0}^{t-\theta}\int_{\{z_n=s\}}G(x-y,t-s)\, \bar g(v_n(y,s))\,d{\mathcal H}^{N-1}(y)ds \\
& \quad=  \int_{\{0<z<t-\theta\}}G(x-y,t-z(y))\,  \bar g(\bar v(y,z(y)))|Dz(y)| \,dy \\
& \quad= \int_{0}^{t-\theta}\int_{\{z=s\}} G(x-y,t-s)\, \bar g(\bar v(y,s)) \,d{\mathcal H}^{N-1}(y)ds \; .
\end{aligned}
$$
Since, as above, 
$$
\left| \bar v(x,t)+\kappa\int_{0}^{t-\theta}\int_{\{z=s\}} G(x-y,t-s) \bar g(\bar v(y,s)) d{\mathcal H}^{N-1}(y)ds \right| \leq C(\bar \rho)\, \theta^{1/2}\; ,
$$
we have proved that $\bar v$ satisfies
$$
\bar v(x,t)=-\kappa\int_{0}^{t}\int_{\{z=s\}} G(x-y,t-s) \bar g(\bar v(y,s)) d{\mathcal H}^{N-1}(y)ds\;,
$$
{\it i.e.}, $\bar v$ is a solution to 
$$
\left\{\begin{array}{rl}
v_t-\Delta v+\kappa\bar g(v){\mathcal H}^{N-1}\lfloor\{u(\cdot,t)=0\}=0 & {\rm in }\; \R^N\times (0,T),\\
v(x,0)=0 & {\rm in }\; \R^N.
\end{array}\right.
$$
The solution of this equation being unique, we have $\bar v=v$, which
proves the convergence of $(v_n)$ to $v$. 
\QED

%%%%%%%%%%%%%%%%%%%%%%%%%%%%%%%%%%%%%%%%%%%%%%
\subsection{Proof of the existence Theorem \ref{main}}

We are now ready to prove Theorem \ref{main}. Throughout the proof, 
$C$ denotes a constant which depends only the data of the problem: 
$N$, $T$, $\kappa$, $\bar g$, $u_0$ and $v_0$. Let us fix some constants $\bar C,R,C_1>0$ to be chosen later and
let ${\mathcal V}={\mathcal V}(\bar C,R,C_1)$ be the set of maps $v:\R^N\times [0,T]\to\R$ such 
that $v$ is measurable, $1/2$-H\"{o}lder continuous in space with 
constant $\bar C$, bounded by a constant $R>\|v_0\|_\infty$ and such that
$$
|v(x,t)-v(y,t)|\leq C_1|x-y|(1+|\log |x-y| |) 
\quad \text{ for all } x,y\in \R^N, \ t\in [0,T].
$$
 Notice that 
${\mathcal V}$ is a closed convex subset of the Banach space
$L^\infty(\R^N\times [0,T]).$

To any $v\in {\mathcal V}$ we associate a map $\tilde v$ 
defined in the following way: let $u$ be the solution to 
$$
\left\{\begin{array}{l}
u_t(x,t)=\bar g(v(x,t))|Du(x,t)|\\
u(x,0)=u_0(x),
\end{array}\right.
$$ 
and let us set $$K(t)=\{u(\cdot,t)\geq0\}, \quad \Gamma(t)=\partial K(t) \quad \text{and} \quad z(x)=\inf\{t\geq 0\; ; \; x\in K(t)\}.$$

Since the velocity $c(x,t):= \bar g(v(x,t))$ satisfies \eqref{BornesVitesse}, \eqref{estiLn} and is 
$1/2$-H\"{o}lder continuous in space 
with constant $\| \bar g' \|_\infty \bar C$,
and since the initial condition enjoys the interior ball property, 
we know from Corollary \ref{cor:ccone} with $\beta=\alpha-1/p=1/2$
that each $K(t)$ has the interior cone property of parameter 
$(\bar \rho,2\bar \rho)$, where $\bar \rho=C_0{\bar C}^{-2}$. Moreover, by \eqref{representation} there exists $M>0$ depending only on the data such that for any $t\in [0;T]$, $K(t)\subset \overline{B}(0,M)$, while \eqref{estiDz} holds thanks to Proposition \ref{zLipschitz}. 

By Lemma \ref{ExStabHeat1} we can therefore define the unique solution $\tilde v$ to 
$$
\left\{\begin{array}{l}
\tilde v_t(x,t)-\Delta \tilde v(x,t)
+\bar g(\tilde v(x,t)){\mathcal H}^{N-1}\lfloor\{u(\cdot,t)=0\}=0\\
\tilde v(x,0)=v_0(x).
\end{array}\right.
$$
From Lemma \ref{ExStabHeat1} we also have, for all   
$x,y\in \R^N,$ $0\leq t\leq t+h\leq T,$
$$
\begin{aligned}
|\tilde v(x,t)|&\leq C(1+|\log(\bar \rho)|)
\leq C(1+|\log(\bar C)|),\\
&\\
|\tilde v(x,t)-\tilde v(y,t)|&\leq 
C(1+|\log(\bar \rho)|)(\bar \rho)^{-1/4} \, |x-y|^{1/2}\leq
C(1+|\log(\bar C)|){\bar C}^{1/2} \, |x-y|^{1/2},\\
&\\
|\tilde v(x,t)-\tilde v(y,t)|&\leq 
\frac{C}{\bar \rho}\, |x-y|(1+|\log(|x-y|)|)\leq C\, {\bar C}^{2} \, |x-y|(1+|\log(|x-y|)|),
\end{aligned}
$$
and
$$
|\tilde v(x,t+h)-\tilde v(x,t)|\leq \frac{C}{\bar \rho} (1+|\log |h| |)\,  h^{1/2}\leq C\, {\bar C}^{2} \, (1+|\log |h| |)\,  h^{1/2}.
$$
So, if we choose  $\bar C$ such that 
$$
C  (1+|\log(\bar C)|){\bar C}^{1/2} \leq \bar C
$$
and then $R$ and $C_1$ such that
$$
R\geq C(1+|\log(\bar C)|) \; {\rm and }\; 
C_1\geq C\, {\bar C}^{2}  \;,
$$
we obtain that $\tilde v \in {\mathcal V}$. Let us now fix $\bar C, R$ 
and $C_1$ as above. Then the map $\Phi,$ which associates to 
$v\in {\mathcal V}$ the map $\tilde v,$ is compact because of the $L^\infty$
and H\"{o}lder bounds on $\tilde v$ recalled above. Since, from Lemma \ref{continuite}, $\Phi$ is also continuous, 
we can complete the proof thanks to Schauder's fixed point theorem.
\QED

%%%%%%%%%%%%%%%%%%%%%%%%%%%%
%%%%%%%%%%%%%%%%%%%%%%%%%%%%

\end{document}